\documentclass[11pt]{amsart}
\pdfoutput=1

\usepackage{eucal}
\usepackage{outlines}
\usepackage{amssymb,amsfonts}
\usepackage{stmaryrd}

\usepackage{autobreak}
\usepackage{enumerate}
\usepackage{mathrsfs}
\usepackage{tikz-cd}
\usepackage{mathtools}
\usepackage{bm}
\usepackage{quiver} 
\usepackage{scalerel}
\usepackage{stackengine,wasysym}
\usepackage{csquotes}
\usepackage{amsmath,calligra,mathrsfs}
\usepackage[style=alphabetic,
minalphanames=1,
maxalphanames=7,
maxbibnames=99]{biblatex}
\usepackage{hyperref}
\usepackage[capitalise]{cleveref}

\newtheorem{thm}{Theorem}[section]

\newtheorem{cor}[thm]{Corollary}
\newtheorem{prop}[thm]{Proposition}
\newtheorem{lem}[thm]{Lemma}

\newtheorem{claim}[thm]{Claim}

\theoremstyle{definition}

\newtheorem{defn}[thm]{Definition}

\newtheorem{exmp}[thm]{Example}

\newtheorem{notn}[thm]{Notation}

\newtheorem{disc}[thm]{Discussion}

\newtheorem{warn}[thm]{Warning}
\newtheorem{remark}[thm]{Remark}
\newtheorem{conv}[thm]{Convention}
\newtheorem{assumption}[thm]{Assumption}

\newtheorem*{mainidea*}{Main Idea}
\newtheorem*{disc*}{Discussion}

\theoremstyle{remark}

\newtheorem{rem}[thm]{Remark}

\newtheorem{notation}[thm]{Notation}

\newcommand{\beqn}{\begin{equation}}
	\newcommand{\eeqn}{\end{equation}}
\newcommand{\bclaim}{\begin{claim}}
	\newcommand{\eclaim}{\end{claim}}
\newcommand{\blem}{\begin{lem}}
	\newcommand{\elem}{\end{lem}}
\newcommand{\bproof}{\begin{proof}}
	\newcommand{\eproof}{\end{proof}}
\newcommand{\bdef}{\begin{defn}}
	\newcommand{\edefn}{\end{defn}}
\newcommand{\bprop}{\begin{prop}}
	\newcommand{\eprop}{\end{prop}}
\newcommand{\bthm}{\begin{thm}}
	\newcommand{\ethm}{\end{thm}}
\newcommand{\brem}{\begin{rem}}
	\newcommand{\erem}{\end{rem}}
\newcommand{\bcor}{\begin{cor}}
	\newcommand{\ecor}{\end{cor}}

\newcommand{\bbC}{\mathbb{C}}

\newcommand{\bbE}{\mathbb{E}}

\newcommand{\bbN}{\mathbb{N}}

\newcommand{\bbR}{\mathbb{R}}

\newcommand{\bbZ}{\mathbb{Z}}

\newcommand{\calA}{\mathcal{A}}

\newcommand{\calC}{\mathcal{C}}
\newcommand{\calD}{\mathcal{D}}

\newcommand{\calM}{\mathcal{M}}

\newcommand{\calO}{\mathcal{O}}
\newcommand{\calP}{\mathcal{P}}

\newcommand{\calS}{\mathcal{S}}
\newcommand{\calT}{\mathcal{T}}
\newcommand{\calU}{\mathcal{U}}
\newcommand{\calV}{\mathcal{V}}

\newcommand{\calX}{\mathcal{X}}

\newcommand{\scrA}{\mathscr{A}}

\newcommand{\scrC}{\mathscr{C}}
\newcommand{\scrD}{\mathscr{D}}

\newcommand{\on}{\operatorname}
\newcommand{\wit}{\widetilde}
\newcommand{\under}{\underline}
\newcommand{\un}{\underline}
\newcommand{\wih}{\widehat}

\newcommand{\C}{\mathbb{C}}
\newcommand{\R}{\mathbb{R}}

\newcommand{\hook}{\hookrightarrow}

\renewcommand{\-}{\textendash}

\mathchardef\mhyphen="2D

\newcommand{\heart}{\heartsuit}

\newcommand{\oblv}{\mathrm{oblv}}

\newcommand{\Cat}{\mathcal{C}\mathrm{at}}

\newcommand{\End}{\mathrm{End}}

\newcommand{\CAlg}{\mathrm{CAlg}}

\newcommand{\Vect}{\mathsf{Vect}}

\newcommand{\Perf}{\on{\mathsf{Perf}}}

\newcommand{\Shv}{\on{Shv}}

\newcommand{\Fun}{\on{Fun}}

\newcommand{\op}{\on{op}}

\renewcommand{\Bar}{\on{Bar}}

\newcommand{\Cond}{\on{Cond}}
\newcommand{\Ani}{\on{Ani}}
\newcommand{\Liq}{\on{Liq}}
\newcommand{\uHom}{\under{\Hom}}
\newcommand{\Ab}{\on{Ab}}

\newcommand{\proet}{\textrm{pro\'{e}t}}

\DeclareMathOperator{\Mod}{Mod}

\DeclareMathOperator{\Hom}{Hom}
\DeclareMathOperator{\Sp}{Sp}

\DeclareMathOperator{\sHom}{\mathscr{H}\text{\kern -3pt {\calligra\large om}}\,}

\newcommand{\ExtProf}{\on{ExtProf}}
\newcommand{\Cont}{\on{Cont}}
\newcommand{\Prof}{\on{Prof}}

\renewcommand{\Pr}{\calP\mathrm{r}}
\renewcommand{\over}{\overline}
\newcommand{\Ban}{\mathbf{Ban}}
\newcommand{\Groth}{\mathsf{Groth}}
\newcommand{\ab}{\mathrm{ab}}
\newcommand{\sep}{\mathrm{sep}}
\newcommand{\lex}{\mathrm{lex}}
\newcommand{\Ch}{\on{Ch}}

\newcommand{\TVect}{\calT\calV\mathrm{ect}}
\newcommand{\lc}{\mathrm{l.c.}}
\newcommand{\comp}{\mathrm{comp}}
\renewcommand{\disc}{\mathrm{disc}}
\newcommand{\Top}{\calT\mathrm{op}}

\newcommand{\dg}{\mathrm{dg}}
\newcommand{\proj}{\mathrm{proj}}
\newcommand{\rmN}{\mathrm{N}}

\newcommand{\Set}{\calS\mathrm{et}}
\newcommand{\cp}{\mathrm{c.p.}}
\newcommand{\ch}{\mathrm{ch}}
\newcommand{\LFun}{\mathrm{LFun}}

\newcommand{\Pro}{\mathsf{Pro}}
\newcommand{\Fin}{\mathcal{F}\mathrm{in}}
\newcommand{\Pyk}{\mathrm{Pyk}}
\newcommand{\PShv}{\mathcal{P}\mathrm{Shv}}
\newcommand{\Ring}{\mathrm{Ring}}
\newcommand{\Grp}{\mathrm{Grp}}
\newcommand{\TRing}{\calT\Ring}
\newcommand{\TGrp}{\calT\Grp}

\newcommand{\LMod}{\mathrm{LMod}}

\newcommand{\Berk}{\mathrm{Berk}}

\DeclareFontFamily{U}{mathx}{\hyphenchar\font45}
\DeclareFontShape{U}{mathx}{m}{n}{
      <5> <6> <7> <8> <9> <10>
      <10.95> <12> <14.4> <17.28> <20.74> <24.88>
      mathx10
      }{}
\DeclareSymbolFont{mathx}{U}{mathx}{m}{n}
\DeclareFontSubstitution{U}{mathx}{m}{n}
\DeclareMathAccent{\widecheck}{0}{mathx}{"71}
\DeclareMathAccent{\wideparen}{0}{mathx}{"75}

\makeatletter
\newcommand*\bigcdot{\mathpalette\bigcdot@{.5}}
\newcommand*\bigcdot@[2]{\mathbin{\vcenter{\hbox{\scalebox{#2}{$\m@th#1\bullet$}}}}}
\makeatother

\usepackage{geometry}
\title{Liquid functional calculus}

\author{Kendric Schefers}

\begin{document}

\begin{abstract}
We develop an elementary formalism of functional
calculus for entire holomorphic functions
in the setting of Clausen and Scholze's $p$-liquid vector spaces.
\end{abstract}

\maketitle

\setcounter{tocdepth}{2}
\tableofcontents

\section{Introduction}

\subsection{Classical functional calculus}
Let $k$ be a field, and let $V$ be a vector space over $k$.
The data of a linear map $T: V \to V$ is equivalent to a choice
of element in the set $\Hom_{\Set}(\ast, \End(V))$, where 
$\End(V)$ denotes the set of linear endomorphisms of $V$.
Of course $\End(V)$ has additional structure: it is canonically 
a $k$-vector space (in fact a $k$-\emph{algebra}). There is an adjunction between
$\Set$ and the category $\Vect_k$\footnotemark
\footnotetext{This notation conflicts with our conventions in
the rest of the paper where $\Vect_k$ will denote the derived $\infty$-category
$k$-vector spaces, or, equivalently, the $\infty$-category of
$k$-module spectra.}
of $k$-vector spaces given by the
forgetful functor from $k$-vector spaces to $\Set$ and the
free $k$-vector space construction on sets, which furnishes an equivalence
\[\Hom_{\Set}(\ast, \End(V)) \simeq \Hom_{\Vect_k}(k[x], \End(V)).\]
This equivalence gives a functional calculus for polynomials---given
a linear operator, you can act by polynomials in this operator.

If $k = \C$, we might consider topological
vector spaces over $\C$. If the topology on such a vector space $V$ is 
induced by a complete metric, such as when $V$ is a Banach
space, we can go further and make sense of things like
everywhere convergent power series, i.e. entire holomorphic
functions $\calO(\C)$, in this operator. In other words,
if $V$ is Banach, the equivalence above lifts as follows,
\beqn
\label{eqn: intro diagram 1}
\begin{tikzcd}
	& {\Hom_{\Vect_{\C}}(\calO(\C), \End(V))} \\
	{\Hom_{\Set}(\ast, \End(V))} & {\Hom_{\Vect_{\C}}(\C[x], \End(V))}
	\arrow[from=1-2, to=2-2]
	\arrow["\simeq", from=2-1, to=2-2]
	\arrow[dashed, from=2-1, to=1-2],
\end{tikzcd}
\eeqn  
where the vertical arrow is restriction
along the inclusion of $\C$-vector spaces
$\C[x] \to \calO(\C)$. This lift defines what we call the
\emph{entire holomorphic functional calculus}.

\subsection{Functional calculus for liquid vector spaces}

\subsubsection{Condensed vector spaces}


Condensed sets are sheaves of sets defined on the site of
profinite sets equipped with the effective epimorphism topology.
Condensed sets (nearly) form a topos, so one may consider all
manner of algebraic objects, such as rings and modules, 
defined internally to this topos. A beautiful feature of this category
is that it contains the category of compactly generated 
topological spaces as a full subcategory. As a consequence,
the category of compactly generated 
topological $\C$-vector spaces embeds fully faithfully into
the category of condensed $\C$-vector spaces, and the latter
is shown by Clausen and Scholze (\cite{Condensed}) 
to be an extremely well-behaved abelian category.

The procedure of the previous section for obtaining a polynomial
functional calculus goes through \emph{mutatis mutandis}
for linear operators in the condensed topos $\Cond$. 
The endomorphisms of a condensed $\C$-vector space
$V$ are a condensed set $\End(V)$ equivalent to the internal mapping object
$\Hom_{\Cond}(\ast, \End(V))$, and by adjunction, we obtain an equivalence,
\[\Hom_{\Cond}(\ast, \End(V)) \simeq \Hom_{\Cond(\C)}(\C[\ast], \End(V)),\]
where $\C[\ast]$ denotes the free condensed $\C$-vector space
on the point and $\Cond(\C)$ denotes the category
of condensed $\C$-vector spaces. If $V$ is in the 
image of the embedding of topological
vector spaces into condensed ones, this equivalence
reduces to the polynomial functional calculus of the previous section.

In the previous section, in order to extend the functional calculus
from polynomials to entire functions, we restricted our attention to
complete topological vector spaces. In the condensed setting,
the role of complete topological vector spaces is played by \emph{$p$-liquid
vector spaces}.

\subsubsection{Liquid vector spaces}
The category of $p$-liquid vector spaces, for
fixed $0 <p \leq 1$, is an abelian category which serves as a well-behaved 
enlargement of the category of Banach spaces over $\C$.
Roughly, it is defined as follows. To each (extremally disconnected\footnotemark)
\footnotetext{See the discussion at the end of \S\labelcref{ssec: Condensed objects}.}
profinite set $S$ is assigned a space of ``measures," denoted
$\calM_{<p}(S)$. One can imagine that an element of $\calM_{<p}(S)$
assigns a weight to each point in $S$. The points of $S$ may themselves
be considered a measure by sending $s \in S$ to the Dirac
measure assigning weight $1$ to $s$ and $0$ to all other
points of $S$; this gives a map $S \to \calM_{<p}(S)$. 
A $p$-liquid vector space $V$ is one for which every 
map of condensed sets $\varphi: S \to V$ extends to a map
$\wit{\varphi}: \calM_{<p}(S) \to V$ of condensed vector spaces.
This extension can be thought of as sending
$\mu \in \calM_{<p}(S)$ to the $S$-indexed sum of the images
of the points $s \in S$ with coefficients given by the weights
of $\mu$---in other words, it can be thought of as sending
a ``measure" $\mu$ to the ``integral" $\int_S \varphi d\mu$.

\subsubsection{The main idea}
Let $V$ be a $p$-liquid $\C$-vector space.
A choice of endomorphism $T: V \to V$ determines
a map of condensed sets $\bbN \to \End(V)$ by the assignment 
$T \mapsto T^n$, where $\End(V)$ is the condensed $\C$-algebra
and the internal mapping object of $\Cond(\C)$.
Assume furthermore that $\End(V)$ is $p$-liquid. 
If we pretend for a moment that $\bbN$ is a profinite set,
the property that $\End(V)$ is $p$-liquid implies
the existence of the dotted line in the following
diagram,
\beqn
\label{eqn: intro diagram 2}
\begin{tikzcd}
	& {\Hom_{\Cond(\C)}(\calM_{<p}(\bbN), \End(V))} \\
	{\Hom_{\Cond}(\ast, \End(V))} & {\Hom_{\Cond}(\bbN,\End(V))}
	\arrow["\simeq", from=1-2, to=2-2]
	\arrow[from=2-1, to=2-2]
	\arrow[dashed, from=2-1, to=1-2],
\end{tikzcd}
\eeqn
because the vertical arrow is an equivalence.
Compare diagram \labelcref{eqn: intro diagram 2} 
to the diagram \labelcref{eqn: intro diagram 1},
in which the classical functional calculus
was defined by the lift indicated by the
dotted line in the latter diagram. The similarity
between these two diagrams suggests the main idea
behind our construction of an extension of the liquid functional
calculus from polynomials to entire holomorphic functions on $\C$.

\begin{mainidea*}
Our strategy for defining a functional calculus for liquid vector spaces
will be to \emph{identify entire functions as a subspace of
$\calM_{<p}(\bbN)$}. Using this approach, for a given operator
$T: \ast \to \End(V)$, the operator $f(T)$
is obtained by evaluating the lift $\wit{T}: \calM_{<p}(\bbN) \to \End(V)$
on the measure $\mu_f \in \calM_{<p}(\bbN)$ associated to the entire
function $f$.
\end{mainidea*}

\brem
Of course, $\bbN$ is not a profinite set, but this problem is
easily solved by considering the space $\bbN \cup \{\infty\}$,
which is (though it is not extremally disconnected).
\erem

\subsection{Results}
We define the measure $\mu_f \in \calM_{<p}(\bbN \cup \{\infty\})$
associated to an entire function $f$ in \cref{def: f as a measure}.
Our definition is actually shown to 
given an element in the space $\calM_{<p}(\bbN \cup \{\infty\})$
in \cref{prop: mu is a measure}.
Using $\mu_f$, we define in
\cref{def: liquid functional calculus} 
the operator $f(T) \in \End(V)$ for 
an endomorphism $T$ of an object $V$ of the derived
$\infty$-category $\calD_{\Ban}(\C, \calM_{<p})$ 
of $p$-liquid vector spaces
with Banach space cohomologies. By restricting
our attention to objects of this category, we guarantee
that its space of endomorphisms is $p$-liquid.

The entire functional calculus defined in this way
turns out to be compatible with the entire classical functional calculus
on Banach spaces. More precisely, in \cref{cor: liquid-classical comparison}
we show that if $V$ is an object in the heart of $\calD_{\Ban}(\C, \calM_{<p})$
(i.e. a Banach space viewed as complex concentrated in
degree $0$) the operator $f(T)$ defined
in \cref{def: liquid functional calculus} coincides with
the image of the bounded linear operator obtained from
the classical functional calculus on Banach spaces
under the embedding of Banach spaces into
$p$-liquid vector spaces. Additionally, we show in
\cref{prop: liquid-classical comparison} that the
induced maps on homology spaces given by $f(T)$
agrees with the maps obtained via the classical functional
calculus on Banach spaces by applying $f$ to 
the induced maps on homology spaces given by $T$.

As a trivial application of our results, we deduce
a functional calculus for perfect complexes of $\C$-vector
spaces in \cref{thm: perfect functional calculus}

\subsubsection{Significance}
The results described above allow one to make sense
of the make sense of things like the exponential of an
automorphism of complexes of Banach spaces, up to
$p$-liquid quasi-isomorphism. The author has no immediate
application in mind for such a device, but examples abound
in mathematical physics where
one considers complexes of infinite-dimensional
vector spaces with non-trivial topologies.

\subsection{An alternative construction}
It was pointed out to the author by P. Scholze that an entire holomorphic functional
calculus is implicitly found in one of the constructions in \cite{ComplexCondensed}.

Namely, in \cite[Lecture 6]{ComplexCondensed} a full subcategory of
$C(\C[T]) = \calD(\Liq_p(\C[T]))$ called $C(\C[T],\C[T])$ is defined
consisting of precisely those objects whose $\C[T]$-actions 
extends uniquely to an action of the entire holomorphic functions in $T$.
One can then wonder whether there are conditions on an object $V \in \calD(\Liq_p(\C))$ 
that ensure that for any endomorphism $T: V \to V$, the corresponding
object of $\calD(\Liq_p(\C[T]))$ will lie in the subcategory 
$C(\C[T],\C[T])$. As we observe in this work, it is enough to assume that all cohomologies of $V$ 
are Banach, which is also recorded (on the abelian level) by \cite[Proposition 5.13]{ComplexCondensed}.

The approach suggested by Scholze is, in fact, more general than our own, 
since it straightforwardly yields a functional calculus for holomorphic 
functions defined on only a subset of the spectrum of an operator, 
as we informally indicate below.

Note that $C(\C[T],\C[T])$ localizes on $\C$ via the map from the locale 
$\calS(\C[T], \C[T])$ to the Berkovich spectrum $\calM^{\Berk}(\C[T]) \simeq \C$. Sections of this sheaf
over an open subset $U \subset \C$ include $p$-liquid $\C[T]$-modules whose 
action of $\C[T]$ extends to an action of the holomorphic functions on $U$. 
By \cite[Proposition 5.14]{ComplexCondensed}, the support of an object
$M \in C(\C[T],\C[T])$ as a section of this sheaf coincides with the 
spectrum of the operator $T: M \to M$ (at least on the abelian level).
Thus an object $V \in \calD(\Liq_p(\C[T]))$ belonging to $C(\C[T],\C[T])$
belongs to $C(\C[T],\C[T])(U)$ if and only if $U$ contains the spectrum
of $T$; that is, the induced action of $\C[T]$ on a complex $V \in \calD(\Liq_p(\C))$ equipped
with an endomorphism extends to an action of holomorphic functions on $U$ precisely when 
$U$ contains the spectrum of $T$.

\section{Notation}

\begin{notn}
\label{notn: Mod}
If $R$ is a classical associative ring, we may view it as
an $\bbE_1$-algebra in $\Sp$. We let $\LMod_{R}$
denote the $\infty$-category of left $R$-module objects
in $\Sp$. This category has a canonical $t$-structure,
with respect to which $\LMod_R^{\heart}$ is the abelian
category of left $R$-modules. Note that
$\LMod_R$ is equivalent to the unbounded 
derived $\infty$-category
of its heart, $\calD(\LMod_R^{\heart})$.

When $R$ is commutative, it has a natural
structure of an $\bbE_{\infty}$-algebra,
and we let $\Mod_{R}$ denote the category of 
$R$-modules in spectra.
\end{notn}

\subsection{Morphisms}

We use $\Prof$ to denote the category
of profinite sets, i.e.\ compact Hausdorff totally disconnected
topological spaces.

We use the notation ``$\uHom_{\scrC}(-,-)$" to
denote the internal hom of a category $\scrC$;
that is, $\uHom_{\scrC}(-,-)$ is always itself an object
of $\scrC$.

On the other hand, we use the notation ``$\Hom_{\scrC}(-,-)$"
to denote the $\infty$-groupoid (resp.\ set) of morphisms in the $\infty$-category
(resp.\ 1-category) $\scrC$.

\subsection{Non-abelian derived categories}
We use $\Ani$ to denote the $\infty$-category of spaces,
the non-abelian derived category of finite sets.

\subsection{Grothendieck universes}
It is inevitable when working with higher
categories that one will encounter ``large"
collections of objects,
i.e. collections that do not form sets.
We adopt the same approach these objects
as Lurie outlines in \S1.2.15 of \cite{HTT},
and which we recall below.

\subsubsection{Strongly inaccessible cardinals}

\begin{assumption}
We assume that for each cardinal $\alpha_0$,
there exists a strongly inaccessible cardinal
$\alpha \geq \alpha_0$.
\end{assumption}

Let $\calU(\alpha)$ denote the collection of all
sets having rank $<\alpha$. Then $\calU(\alpha)$
is a Grothendieck universe; it satisfies all of the usual axioms
of set theory.

\bdef
A mathematical object is $\alpha$-\emph{small} if
it belongs to $\calU(\alpha)$. It is \emph{essentially
$\alpha$-small} if it is equivalent (in whatever relevant sense)
to an $\alpha$-small object. We let $\Set_{\alpha}$ denote the
category of $\alpha$-small sets.
\edefn

Outside of foundational work, mathematics takes
place in a fixed Grothendieck universe, so
when an author like Lurie writes ``small,"
he has implicitly chosen a Grothendieck
universe.

\subsection{Glossary of categories}

This work involves many different categories
whose definitions and differences are often
subtle. This work also contains many
results which may be of independent
interest to the reader. To facilitate the reader
trying to read this work in a piecemeal fashion,
we furnish a glossary of categories below. 
We have tried to group the categories by
theme. In what follows, $\kappa$ is an
uncountable strong limit cardinal,
and $\alpha > \kappa$ is a strongly
inaccessible cardinal.

\subsubsection{Glossary of categories}
${}$ \\

\begin{outline}
\0 \noindent Categories of purely algebraic objects.

\1 $\Ring$ is the ordinary category of rings.
	\2$\Ring_{\alpha}$ is the full subcategory of $\Ring$
	on $\alpha$-small rings.

\1 $\Grp$ is the ordinary category of groups.
	\2 $\Grp_{\alpha}$ is the full subcategory of $\Grp$
	on $\alpha$-small groups.

\1 $\Ab$ is the ordinary category of abelian groups.
	\2 $\Ab_{\alpha}$ is the full subcategory of $\Ab$ on
	$\alpha$-small abelian groups.

\1 $\calD_{\geq 0}(\Ab)$ is the $\infty$-category of non-negative
	homologically graded chain complexes in $\Ab$. It is equivalent
	to the animation of $\Ab$, $\Ani(\Ab)$.
	\2 $\calD(\Ab)$ is the stabilization of $\calD_{\geq 0}(\Ab)$,
	also known as the unbounded derived $\infty$-category of $\Ab$.

\1 $\Vect_{\C}$ is the unbounded derived $\infty$-category of $\C$-vector
spaces, also know as the $\infty$-category of chain complexes of
	$\C$-vector spaces with quasi-isomorphisms inverted.
	\2 $\Vect_{\C}^{f.d.}$ is ordinary category of finite dimensional
	$\C$-vector spaces. They are the compact projective generators
	of $\Vect_{\C}$.

\0 Categories of categories

\1 $\wih{\Cat}_{\infty}$ is the category of \emph{all}
(i.e. not necessarily small) $\infty$-categories.
	\2 $\Cat_{\infty}$ is the category of small $\infty$-categories.

\0 Categories of topological objects.

\1 $\Top$ is the ordinary category of topological spaces,
whose morphisms are continuous functions. 
	\2 $\Top_{\kappa}$ is the full subcategory of $\Top$ on
	$\kappa$-compactly generated topological spaces.

\1 $\TRing$ is the ordinary category of topological rings,
whose morphisms are continuous ring homomorphisms.
	\2 $\TRing_{\kappa}$ is the full subcategory of $\TRing$
	on objects whose underlying topological space is $\kappa$-compactly
	generated.

\1 $\TGrp$ is the ordinary category of topological groups,
whose morphisms are continuous group homomorphisms.
	\2 $\TGrp_{\kappa}$ is the full subcategory of $\TGrp$ on
	objects whose underlying topological space is $\kappa$-compactly
	generated.

\1 $\TVect$ is the ordinary category of topological
$\C$-vector spaces, whose morphisms are continuous
linear maps.
	\2 $\TVect_{\lc}^{\comp}$ is the full subcategory of $\TVect$
	on locally compact topological vector spaces whose topology
	is induced by a complete metric.

\1 $\Ban$ is the ordinary category of complex Banach spaces,
	whose morphisms are bounded linear operators.

\0 Categories coming from condensed mathematics.

\1 $\Cond$ is the category of $\kappa$-condensed sets.
	\2 $\Cond(\Ring)$ is the category of $\kappa$-condensed rings.
	\2 $\Cond(\Grp)$ is the category of $\kappa$-condensed groups.
	\2 $\Cond(\Ab)$ is the category of $\kappa$-condensed abelian groups.

\1 $\Cond(\scrC)$ is the category of $\kappa$-condensed
objects in the $\infty$-category $\scrC$.

\1 $\Cond(\calA)$ is the category of modules in $\Cond(\Ab)$ 
over the condensed commutative ring $\calA$.
	\2 $\heart(\calA, \calM)$ is the full subcategory of $\Cond(\calA)$
	determined by the analytic ring structure $(\calA, \calM)$, defined in \cref{def: discrete
	analytic modules}.

\1 $\calD_{\geq 0}(\calA)$ is the $\infty$-category of modules in $\Cond(\calD_{\geq 0}(\Ab))$
over the condensed animated ring $\calA$. 
	\2 $\calD_{\geq 0}(\calA, \calM)$ is the full subcategory of $\calD_{\geq 0}(\calA)$
	determined by the analytic animated ring structure $(\calA, \calM)$, defined in \cref{def: analytic modules}.

\1 $\calD(\calA)$ is the stable $\infty$-category given by 
the stabilization of $\calD_{\geq 0}(\calA)$.
	\2 $\calD(\calA, \calM)$ is the stabilization of $\calD_{\geq 0}(\calA, \calM)$.

\1 $\Liq_p(\C)$ is alternative notation for $\heart(\C, \calM_{<p})$, where
$(\C, \calM_{<p})$ denotes the $p$-liquid analytic ring structure on $\un{\C}$.

\1 $\calD(\C, \calM_{<p})$ is the category $\calD(\calA, \calM)$ listed
above for the analytic animated ring $(\C, \calM_{<p})$.
	\2 $\calD_{\Ban}(\C, \calM_{<p})$ is the full subcategory of $\calD(\C, \calM_{<p})$
	on objects whose homology spaces lie in the essential image of the embedding
	$\Ban \hook \Liq_p(\C)$.

\end{outline}



\section{Condensed mathematics}

\subsection{Condensed sets}
Roughly speaking, a condensed set is a
sheaf of sets on the pro-\'etale site of
a geometric point, denoted 
$\ast_{\proet}$. Explicitly, 
$\ast_{\proet}$ is the site given by the category
of profinite sets, denoted $\Prof$,
with open covers given by
finite families of jointly surjective maps.
For the sake of clarity, we recall the
precise definition.

\bdef
Let $\calU(\alpha)$ be a fixed Grothendieck
universe. For any $\alpha$-small, ordinary category $\scrC$,
the category $\Pro^{\alpha}(\scrC)$ is defined as the
full subcategory of $\Fun(\scrC, \Set_{\alpha})^{\op}$
on those functors which are limits of
cofiltered diagrams of functors representable under the
Yoneda embedding $\scrC \hook \Fun(\scrC, \Set_{\alpha})^{\op}$.
We set $\Prof^{\alpha} := \Pro^{\alpha}(\Fin)$.
\edefn

\brem
As a category, $\Prof^{\alpha}$ is equivalent to
the category of compact Hausdorff totally 
disconnected topological spaces whose underlying
sets are $\alpha$-small.
\erem


There is a problem, however, naively defining
condensed sets as sheaves on this site. Given any
choice of Grothendieck universe $\calU(\alpha)$, 
the category $\Prof^{\alpha}$
is large, so it is not a good idea to
work with the category of sheaves on it
since such a category would not be a topos.
Clausen-Scholze circumvent this problem
by working with a modification
of $\Prof^{\alpha}$ obtained as follows. 
Choose an uncountable strong limit cardinal $\kappa < \alpha$,
and instead consider the category of $\kappa$-small
profinite sets, denoted $\Prof_{\kappa}$,
rather than $\Prof^{\alpha}$.

\bdef[{\cite[Definition 2.1]{Condensed}}]
The site $\ast_{\kappa\-\proet}$ is the site of $\kappa$-small 
profinite sets $S$ with covers given
by finite families of jointly surjective maps.
\edefn

\brem
\label{rem: small site}
Note that $\Prof_{\kappa}$ is an $\alpha$-small
category; this follows from the definition of
strongly inaccessible cardinal.
\erem

Clausen-Scholze define the category of
$\kappa$-condensed sets as the category
of $\Set_{\alpha}$-valued (resp.\ $\Grp_{\alpha}$-valued,
$\Ring_{\alpha}$-valued, $\Ab_{\alpha}$-valued) sheaves on $\ast_{\kappa\-\proet}$,
which they denote by $\Cond_{\kappa}$
(resp.\ $\Cond_{\kappa}(\Grp)$, $\Cond_{\kappa}(\Ring)$, $\Cond_{\kappa}(\Ab)$).
For any two choices of uncountable strong
limit cardinal $\kappa' > \kappa$, they show
that there is a fully faithful embedding
of $\kappa$-condensed sets into $\kappa'$-condensed
sets (\cite[Proposition 2.9]{Condensed}). They then
define the category of condensed sets to be the
filtered colimit,
\[\varinjlim_{\kappa} \Cond_{\kappa}\]
(resp.\ $\varinjlim_{\kappa} \Cond_{\kappa}(\Grp)$, $\varinjlim_{\kappa} \Cond_{\kappa}(\Ring)$,
$\varinjlim_{\kappa} \Cond_{\kappa}(\Ab)$),
taken in some suitable category of categories.

\begin{conv}
In the sequel,
we fix a Grothendieck universe $\calU(\alpha)$
by fixing a strongly inaccessible
cardinal $\alpha$. Having fixed a
Grothendieck universe, we omit
any reference to $\alpha$ in our notation.
For example, $\Set_{\alpha}$, $\Pro^{\alpha}$,
and $\Prof^{\alpha}$, respectively, will be denoted
by $\Set$, $\Pro$, and $\Prof$, respectively.
\end{conv}

\subsubsection{Condensed sets as replacements for topological spaces}
Clausen and Scholze defined 
condensed sets as replacements for
topological spaces that have better behavior
when considered with algebraic structures.
Strong evidence for their suitability as a replacement
is given by the following proposition of Clausen-Scholze.

\bprop[{\cite[Proposition 1.7]{Condensed}}]
\label{prop: Proposition 1.7}
Let $X$ be a topological space, and denoted by
$\un{X}$ the condensed set given by the assignment
\[S \mapsto \Cont(S, X),\]
for $S \in \Prof_{\kappa}$, where
$\Cont(S,X)$ denotes the set of continuous
functions $S \to X$. Then $\un{(-)}: \Top_{\kappa}
\to \Cond_{\kappa}$ is a fully faithful
functor from the category of $\kappa$-compactly
generated topological spaces to $\kappa$-condensed
sets.
\eprop

\brem
\label{rem: TRing embeds into Cond(Ring)}
The functor $\un{(-)}$ also induces fully faithful
embeddings of $\TGrp_{\kappa}$ and $\TRing_{\kappa}$
into $\Cond_{\kappa}(\Grp)$ and $\Cond_{\kappa}(\Ring)$, 
respectively.
\erem

\brem
If $R$ is a topological ring with the discrete
topology, $\un{R}$ is a constant sheaf
on $\Prof_{\kappa}$. Note that the image,
for example, of the topological field $\bbR$
under $\un{(-)}$,
however, is \emph{not} a constant sheaf.
In particular, modules in
$\Cond_{\kappa}(\Ab)$ over each of $\un{\bbR}$ and
$\un{\bbR_{\disc}}$ are
different.
\erem

Unlike the category of topological abelian groups,
the category of condensed abelian groups forms
an abelian category with extremely nice properties.

\bthm[{\cite[Theorem 2.2]{Condensed}}]
The category of $\kappa$-condensed abelian groups 
is an abelian category which satisfies Grothendieck’s 
axioms (AB3), (AB4), (AB5), (AB6), (AB3*) and (AB4*), 
to wit: all limits (AB3*) and colimits (AB3) exist, 
arbitrary products (AB4*), arbitrary direct sums (AB4) 
and filtered colimits (AB5) are exact, and (AB6) for 
any index set $J$ and filtered categories $I_j$, $j \in J$,
with functors $i \mapsto M_i$ from $I_j$ to $\kappa$-condensed 
abelian groups, the natural map
\[\varinjlim_{(i_j \in I_j)_j} \prod_{j \in J} M_{i_j} \to \prod_{j \in J} \varinjlim_{i_j \in I_j} M_{i_j}\]
is an isomorphism. Moreover, the category of 
$\kappa$-condensed abelian groups is generated 
by compact projective objects.
\ethm

\brem
A variant of the above theorem holds
for the category of \emph{all} condensed
abelian groups, $\varinjlim_{\kappa} \Cond_{\kappa}(\Ab)$.
\erem

\subsection{Condensed objects}
\label{ssec: Condensed objects}
As a means of proving the various properties
of $\Cond$, Clausen-Scholze take advantage
of the fact that condensed sets are determined
by their values on very special types of profinite sets,
called \emph{extremally disconnected} profinite sets.

\bprop[{\cite[Proposition 2.7]{Condensed}}]
\label{prop: Proposition 2.7}
Consider the site of $\kappa$-small extremally 
disconnected profinite sets, denoted $\ExtProf_{\kappa}$, 
with covers given by finite families of jointly surjective maps. 
Its category of sheaves is equivalent to $\kappa$-condensed
sets via restriction from profinite sets.
\eprop

Motivated by the proposition,
Clausen-Scholze define $\kappa$-condensed objects
of a category $\scrC$, more generally, as follows.

\bdef[{\cite[Definition 11.7]{AnalyticCondensed}}]
Let $\scrC$ be an $\infty$-category that admits
all small colimits.
The category of $\kappa$-condensed objects
of $\scrC$, denoted $\Cond_{\kappa}(\scrC)$, is the
category of contravariant sheaves from $\ExtProf_{\kappa}$
to $\scrC$ that take finite coproducts to products.
\edefn

\brem
If $\scrC$, in addition, admits all small limits, then
$\Cond_{\kappa}(\scrC)$ can be identified
with the category of $\scrC$-valued sheaves
on $\ExtProf_{\kappa}$.
\erem

As above, for any two choices of cardinals, $\kappa' > \kappa$,
there is a fully faithful functor,
\[\Cond_{\kappa}(\scrC) \to \Cond_{\kappa'}(\scrC),\]
given as the left adjoint to the forgetful
functor from $\Cond_{\kappa'}(\scrC)$ to $\Cond_{\kappa}(\scrC)$.
We may likewise form the category of condensed
objects in $\scrC$ as the filtered colimit,
\[\varinjlim_{\kappa} \Cond_{\kappa}(\scrC),\]
taken in, $\wih{\Cat}_{\infty}$, the 
$\infty$-category of all (e.g.\ not necessarily small) 
$\infty$-categories.

\subsubsection{Pyknotic objects}
While Clausen and Scholze were developing
their theory of condensed sets, Barwick and Haine
had been studying essentially the same notion, which
they call \emph{pyknotic sets}. The difference
between pyknotic sets and condensed sets is purely
set-theoretic in nature. We refer the reader to
\S0.3 of \cite{pyknotic} a discussion on the
differences between the two theories.

For our purposes, it will be useful to deploy some
results from theory of pyknotic objects in the current
setting of condensed objects. In order to do so, we
briefly recall the notion of pyknotic set. 
In the notation of the previous section:
assume the existence of a smallest
strongly inaccessible cardinal $\alpha' > \alpha$.
The category $\Prof_{\alpha}$ is small in the universe 
$\calU(\alpha')$, so sheaves on it form a topos in this
larger universe.

\bdef
The category of pyknotic sets is defined
to be the category of $\Set_{\alpha'}$-valued
sheaves on $\Prof^{\alpha}$. More generally,
given an $\infty$-category $\scrC$, the
category of pyknotic objects of $\scrC$,
denoted $\Pyk_{\alpha}(\scrC)$, is
subcategory of contravariant functors
from $\ExtProf^{\alpha}$ to $\scrC$ that
take finite coproducts to products.
\edefn

\brem
As with condensed objects, it $\scrC$ admits
all small limits, then
\[\Pyk_{\alpha}(\scrC) = \Shv_{\scrC}(\ExtProf^{\alpha}).\]
\erem

Though we have fixed $\alpha$ and $\alpha'$
both to be strongly inaccessible cardinals, the
results in the pyknotic literature hold
equally well if we take $\alpha$ to be an uncountable
strong limit cardinal and $\alpha' > \alpha$ a strongly
inaccessible cardinal. In other words, we may take
$\alpha$ to be $\kappa$ of the previous section,
and $\alpha'$ to be $\alpha$ of the previous section. With these
choices, it is clear that
\[\Cond_{\kappa}(\scrC) = \Pyk_{\kappa}(\scrC).\]
As such, we freely use the results of
\cite{pyknotic} below when working with
$\kappa$-condensed objects.

\subsubsection{Condensed objects as sheaves on $\Prof_{\kappa}$}

We would like to prove an analogue of
\cite[Proposition 2.7]{Condensed} for
sheaves on $\Prof_{\kappa}$ with values
in categories other than $\Set$.
As observed in \cref{rem: small site},
$\ExtProf_{\kappa}$ and $\Prof_{\kappa}$ 
are small sites, so each determines
an $\infty$-topos.
By abuse of notation, we let $\Prof_{\kappa}$
and $\ExtProf_{\kappa}$ denote the $\infty$-topoi
determined by each site.

\bprop
Let $\wih{\Prof}_{\kappa}$ denote the
hypercompletion of $\Prof_{\kappa}$.
Then $\wih{\Prof}_{\kappa}$ and $\ExtProf_{\kappa}$
are equivalent $\infty$-topoi.
\eprop

\bproof
This is \cite[Warning 2.2.2]{pyknotic} and
\cite[Corollary 2.4.4]{pyknotic}.
\eproof

\bcor
\label{cor: sheaves on prof vs extprof}
Let $\scrC$ be a presentable $\infty$-category.
Restriction to extremally disconnected
profinite sets determines an equivalence
of categories,
\[\Shv_{\scrC}(\wih{\Prof}_{\kappa}) \xrightarrow{\simeq} \Shv_{\scrC}(\ExtProf_{\kappa}).\]
In particular, if $\scrC$ additionally
admits all small limits, then
\[\Shv_{\scrC}(\wih{\Prof}_{\kappa}) \simeq \Cond_{\kappa}(\scrC).\]
\ecor

\bproof
Not that the category of $\scrC$-valued sheaves
on \emph{any} $\infty$-topos $\calX$ is given by
the Lurie tensor product of presentable $\infty$-categories,
$\calX \otimes \scrC$. Since, $\wih{\Prof}_{\kappa} \xrightarrow{\simeq} \ExtProf_{\kappa}$
via restriction by the above proposition, it follows
that 
\[\Shv_{\scrC}(\wih{\Prof}_{\kappa}) \simeq \wih{\Prof}_{\kappa} 
	\otimes \scrC \xrightarrow{\simeq} \ExtProf_{\kappa}
		\otimes \scrC \simeq \Shv_{\scrC}(\ExtProf_{\kappa}),\]
as desired.
\eproof

Note that since $\wih{\Prof}_{\kappa} \to \Prof_{\kappa}$
is fully faithful, the functor $\Shv_{\scrC}(\wih{\Prof}_{\kappa})
\to \Shv_{\scrC}(\Prof_{\kappa})$ is also fully
faithful.
The significance of \cref{cor: sheaves on
prof vs extprof} is that it allows us to take
sections of condensed objects of $\scrC$
on arbitrary $\kappa$-small profinite sets
though condensed objects are a priori defined
only on extremally disconnected profinite sets,
as shown in the following lemma.

\blem
Suppose that $\scrC$ is a complete presentable
$\infty$-category. Then there is a fully faithful
embedding,
\[\Cond_{\kappa}(\scrC) \hook \Shv_{\scrC}(\Prof_{\kappa}).\]
\elem

\bproof
By \cref{cor: sheaves on prof vs extprof}, it suffices
to show that $\Shv_{\scrC}(\wih{\Prof}_{\kappa})$
embeds fully faithfully into $\Shv_{\scrC}(\wih{\Prof}_{\kappa})$.
First note that the inclusion $i_*: \wih{\Prof}_{\kappa} \subset \Prof_{\kappa}$
is a fully faithful geometric morphism 
(this is true for the inclusion of the hypercomplete objects
of any $\infty$-topos). As such, it admits a left exact left
adjoint $i^*: \Prof_{\kappa} \to \wih{\Prof}_{\kappa}$.
Pointwise composition with $(i^*)^{\op}$ determines a functor,
\[\Fun({\wih{\Prof}_{\kappa}}^{\op}, \scrC) \xrightarrow{i_*} \Fun(\Prof_{\kappa}^{\op}, \scrC)\]
which is fully faithful because $i^*$ is fully faithful.
Moreover, since $i^*$ is left adjoint, $(i^*)^{\op}$ preserves 
small limits, meaning $i_*$ restricts to a functor,
\[\Shv_{\scrC}(\wih{\Prof}_{\kappa}) \xhookrightarrow{i_*} \Shv_{\scrC}(\Prof_{\kappa}),\]
as desired.
\eproof


\subsection{Important convention}
Throughout the remainder of this paper,
we fix an uncountable strong limit cardinal
$\kappa$ and work with objects in 
$\Cond_{\kappa}(\scrC)$. In doing
so, we will omit $\kappa$ from
all of our notation. In particular,
we will use $\Cond(\scrC)$, $\Prof$,
and $\ExtProf$ to denote
$\Cond_{\kappa}(\scrC)$,
$\Prof_{\kappa}$, and $\ExtProf_{\kappa}$,
respectively. When we say ``condensed"
we mean ``$\kappa$-condensed."

\begin{warn}
This is a significant, and potentially confusing,
departure from the conventions of Clausen-Scholze
and the rest of the condensed mathematics literature; 
$\Cond(\scrC)$ universally denotes the colimit
$\varinjlim_{\kappa} \Cond_{\kappa}(\scrC)$.
Thus, for the sake of clarity, we reiterate: $\Cond(\scrC)$
for us denotes the category of $\kappa$-condensed
objects of $\scrC$, for some fixed $\kappa$.
\end{warn}

\subsection{Miscellanea}

We state and prove below a handful of
results that will be useful for our purposes.

\subsubsection{$\Top$ embeds fully faithfully into condensed anima}
In the next section, we would like to view the topological
fields $\bbR$ and $\bbC$ not only as condensed rings under
the embedding of \cref{rem: TRing embeds into Cond(Ring)},
but also as condensed \emph{animated} (meaning simplicial)
rings.

\blem
\label{lem: Top embeds into Cond(Ani)}
There is a fully faithful functor,
\[\Top \to \Cond(\Ani).\]
\elem

\bproof
Given any cocomplete, compactly generated category $\scrC$,
there is an inclusion $\scrC \to \Ani(\scrC)$.
Letting $\scrC = \Cond$, we obtain a fully faithful
functor $\Cond \to \Ani(\Cond)$. Using 
\cite[Proposition 11.8]{AnalyticCondensed}, this gives
a fully faithful functor
$\Cond \to \Cond(\Ani)$.
The desired functor is then obtained by composing
this functor with the embedding $\Top \hook \Cond$
from \cref{prop: Proposition 1.7}.
\eproof

\brem
Obviously, similar results hold for $\TGrp$ and $\TRing$.
\erem

\subsubsection{Stabilization commutes with taking condensed objects}

\blem
\label{lem: sp commutes with cond}
Let $\scrC$ be a presentable $\infty$-category.
Then there is a natural equivalence of $\infty$-categories,
\beqn
\label{eqn: sp commutes with cond}
\Sp(\Cond(\scrC))  \simeq \Cond(\Sp(\scrC)).
\eeqn
\elem

\bproof

Recall that, for any presentable
$\infty$-category $\scrD$, the category 
of $\scrD$-valued sheaves on any 
$\infty$-topos $\calX$
is given by the Lurie tensor product\footnotemark,
\footnotetext{See \cite[\S4.8]{HA} for the definition
of the Lurie tensor product of $\infty$-categories.}
\[\Shv_{\scrC}(\calX) \simeq \scrC \otimes \calX\]
(see \cite[Remark I.1.3.1.6]{SAG}). 
In particular, $\Cond(\Ani)$ is the $\infty$-topos
of sheaves on the small site $\{\ast\}_{\proet_{\kappa}}$, so
$\Shv_{\Sp(\scrC)}(\Cond(\Ani)) \simeq \Sp(\scrC) \otimes \Cond(\Ani)$.
On the other hand, the stabilization of
a presentable $\infty$-category $\scrC$
admits a characterization as the Lurie
tensor product of $\scrC$ with $\Sp$\footnotemark,
\[\Sp(\scrC) \simeq \Sp \otimes \scrC.\]
\footnotetext{See \cite[Example 4.8.1.23]{HA} for a proof of
this fact.} 
Together, these two facts reduce
\labelcref{eqn: sp commutes with cond}
to the claim that,
\[\Sp \otimes (\Cond(\Ani) \otimes \scrC) \simeq \Cond(\Ani) \otimes (\Sp \otimes \scrC).\]
But this is clear from the fact that
the $\infty$-category of presentable
$\infty$-categories, $\Pr^L$, is symmetric monoidal
under the Lurie tensor product.
\eproof

\subsubsection{Stabilization commutes with taking module objects}

\blem
\label{lem: mod commutes with sp}
Suppose that $\calX$ is an $\infty$-topos,
and let $\calO \in \calX$ be
a grouplike commutative algebra object.
Then there is an equivalence of stable
$\infty$-categories,
\[\Sp(\Mod_{\calO}(\calX)) \simeq \Mod_{\calO_{\Sigma}}(\Sp(\calX)),\]
where $\calO_{\Sigma}$ denotes the image
of $\calO$ in $\Sp(\calX)$.
\elem

\bproof
We note that, by \cite[Remark I.1.3.5.1]{SAG},
there is a canonical equivalence,
\[\Shv_{\CAlg(\Sp)}(\calX) \simeq \CAlg(\Shv_{\Sp}(\calX)).\]
The reasoning in that remark applies
to prove to a similar statement for $\CAlg(\calS)$.
Namely the forgetful functor
$\CAlg(\calS) \to \calS$
is conservative and preserves small
limits by \cite[Lemma 3.2.2.6 and Corollary 3.2.2.5]{HA}.
It follows that we have a canonical equivalence,
\[\Shv_{\CAlg(\calS)}(\calX) \simeq \CAlg(\calX).\]

We now note that the pointwise application
of the suspension functor,
$\Sigma^{\infty}: \calS \to \Sp$, induces
a functor, $F_{\Sigma^{\infty}}: 
\PShv_{\calS}(\calX) \to \PShv_{\Sp}(\calX)$,
which descends to a functor of sheaves,
\[\calX \to \Shv_{\Sp}(\calX),\]
given by $L \circ F_{\Sigma^{\infty}}$,
where $L$ is the sheafification functor,
$\PShv_{\Sp}(\calX) \to \Shv_{\Sp}(\calX)$.
This functor is symmetric monoidal with
respect to the Cartesian monoidal structure
on $\calX$ and the smash product monoidal
structure on $\Shv_{\Sp}(\calX)$, so
it induces a functor,
\[\CAlg(\calX) \xrightarrow{\Sigma} \CAlg(\Shv_{\Sp}(\calX)).\]
Let $\calO_{\Sigma}$
denote the image of $\calO \in \CAlg(\calX)$ 
under the above functor. Now we have the
following chain of equivalences:
\begin{align*}
\Sp(\Mod_{\calO}(\calX)) 	&\simeq \Sp(\calX_{\Bar(\calO)/}) \\
					&\simeq \Sp \otimes (\calX_{\Bar(\calO)/}) \\
					&\simeq \Sp(\calX)_{\Sigma(\Bar(\calO))} \\
					&\simeq \Sp(\calX)_{\Bar(\calO_{\Sigma})} \\
					&\simeq \Mod_{\calO}(\Sp(\calX)),
\end{align*}
where we have used the equivalence $\Mod_{\calO}(\calX) \simeq \calX_{\Bar(\calO)/}$
(see \cite[Remark 5.2.6.28]{HA}); the fact that $\Bar(\calO)$ is again
an object of $\CAlg(\calX)$ under the equivalence given by
the forgetful functor $\CAlg(\calX) \xrightarrow{\simeq} \calX$;
and the fact that $\Sigma$ commutes with geometric realizations
as a left adjoint, so $\Bar(\calO_{\Sigma}) \simeq \Sigma(\Bar(\calO))$.
\eproof

\section{Analytic and analytic animated rings}

\subsection{Analytic rings}
Let us recall the notion of analytic ring
found in \cite[Lecture VII]{Condensed}.

\bdef
\label{def: analytic ring}
A pre-analytic ring $(\calA, \calM)$ is a condensed 
ring $\calA$ together with a
functor,
\[\calM[-]: \ExtProf \to \Mod_{\calA}(\Cond(\Ab))\]
taking finite disjoint unions to products, and 
a natural transformation of functors,
$\ExtProf \to \Mod_{\calA}(\Cond(\Ab))$,
\[\calA[S] \to \calM[S].\]

An analytic ring is a pre-analytic ring
$(\calA, \calM)$ such that for any complex,
\[C: \cdots C_i \to \cdots \to C_1 \to C_0 \to 0,\]
of $\calA$-modules
in condensed abelian groups, such that all $C_i$
are direct sums of objects of the form $\calM[T]$ for
varying extremally disconnected $T$, the map
\[R\uHom_{\calA}(\calM[S], C) \to R\uHom_{\calA}(\calA[S], C)\]
of complexes of condensed abelian groups
is a quasi-isomorphism for all extremally disconnected sets $S$. 
\edefn

\begin{remark}
Heuristically, a pre-analytic ring is supposed to be
a ring equipped with a notion
in condensed abelian groups
of ``free" module over that ring, given
by the functor in the definition.
The condition in the definition of analytic ring
specifies that this notion of ``free" module should
be well-behaved: maps into particular kinds of
$\calA$-modules cannot distinguish between
the ``free" modules $\calM[-]$, and free $\calA$-modules
in the category of condensed abelian groups.
\end{remark}

\subsection{Derived category of condensed modules}
There is a related notion of analytic animated ring,
in which both the ring $\calA$ and its space of measures
$\calM$ are allowed to be objects in 
condensed anima. 

\begin{notation}
Let $\calA$ be an animated condensed
ring\footnotemark. Let $\calD_{\geq 0}(\calA)$ denote 
the prestable $\infty$-category,
$\Mod_{\calA}(\Cond(\calD_{\geq 0}(\Ab)))$,
the category of $\calA$-modules in
animated condensed abelian groups.
\footnotetext{By \cite[Lemma 11.8]{AnalyticCondensed},
we may permute the adjectives ``animated" and ``condensed"
with impunity.}

We denote the stabilization of $\calD_{\geq 0}(\calA)$ by
$\calD(\calA) := \Sp(\calD_{\geq 0}(\calA))$. 
It is a stable $\infty$-category with a natural
t-structure whose connective part is $\calD_{\geq 0}(\calA)$.
We observe that $\calD(\calA)$ admits
a forgetful functor $\oblv_{\Sp}: \calD(\calA) \to \Cond(\Sp)$
which we call taking \emph{the
underlying condensed spectrum.}
Indeed, we have the following
chain of equivalences,
\begin{align*}
\calD(\calA) 			&:= \Sp(\Mod_{\calA}(\Cond(\calD_{\geq 0}(\Ab)))) \\
					&\simeq \Mod_{\calA}(\Sp(\Cond(\calD_{\geq 0}(\Ab)))) \tag{\cref{lem: mod commutes with sp}} \\
					&\simeq \Mod_{\calA}(\Cond(\calD(\Ab))) \tag{\cref{lem: sp commutes with cond}} \\
					&\simeq \Mod_{\calA}(\Cond(\Sp)).
\end{align*}
Under this equivalence, $\oblv_{\Sp}$ is
simply the functor on $\Mod_{\calA}(\Cond(\Sp))$
of forgetting the $\calA$-module structure.
\end{notation}

When $\calA$ is commutative, 
both $\calD_{\geq 0}(\calA)$ and $\calD(\calA)$
are canonically symmetric monoidal categories. whose
symmetric monoidal structure is
induced by the commutative algebra $\calA$
(see \cite[Theorem 3.3.3.9]{HA}).

\blem
\label{lem: internal hom}
The canonical symmetric
monoidal structure on $\calD(\calA)$
(resp. $\calD_{\geq 0}(\calA)$), induced by the
commutative algebra $\calA$, is closed. 
That is, the symmetric monoidal product admits a
right adjoint, which we denote by 
$\uHom_{\calD(\calA)}(-,-)$
(resp. $\uHom_{\calA}(-,-)$),
called the \emph{internal hom}.
\elem

\bproof
We note that $\calA$ is a sheaf of $\bbE_{\infty}$-rings
on the $\infty$-topos $\ExtProf$, 
and $\calD(\calA)$ is the category of modules\footnotemark 
over $\calA$ in $\Shv_{\Sp}(\ExtProf)$.
\footnotetext{$\Mod_{\calA}$ in the notation
of Lurie's \cite{SAG}.} As mentioned,
this category has a canonical symmetric monoidal structure,
whose product we denote by $\otimes_{\calA}$.
By \cite[Proposition I.2.1.0.3]{SAG}, $\otimes_{\calA}$
preserves small colimits in both variables, so by
the Adjoint Functor Theorem, it admits a right adjoint,
which is $\uHom_{\calD(\calA)}(-,-)$.

Recall that $\calD_{\geq 0}(\calA)$ is 
tautologically the connective part of the t-structure on  
$\calD(\calA)$, and that $\calA$ is a connective sheaf
of $\bbE_{\infty}$-rings on $\ExtProf$. By
\cite[Proposition I.2.1.1.1(b)]{SAG}, $\calD_{\geq 0}(\calA)$
is closed under the tensor product $\otimes_{\calA}$ and contains the
unit object, so $\calD_{\geq 0}(\calA)$ is symmetric
monoidal under $\otimes_{\calA}$. Moreover,
$\calD_{\geq 0}(\calA)$ is closed under small colimits,
and reflects all colimits in $\calD(\calA)$, so
$\otimes_{\calA}$ preserves small colimits
in $\calD_{\geq 0}(\calA)$, as well. Thus,
$\otimes_{\calA}$ admits a right adjoint
in $\calD_{\geq 0}(\calA)$ called
$\uHom_{\calA}(-,-)$, which, by
uniqueness of the right adjoint, must be
the restriction of $\uHom_{\calD(\calA)}(-,-)$.
\eproof

\brem
A similar argument as presented in
the proof of \cref{lem: internal hom}
proves the existence of an internal
hom for $\Mod_{\calA}(\Cond(\Ab))$
for $\calA$ a discrete condensed 
commutative ring.
\erem

\subsection{Analytic animated rings}

\bdef
\label{def: analytic animated ring}
A pre-analytic animated ring
$(\calA, \calM)$ is an animated condensed
ring $\calA$ together with a functor,
\[\calM[-]: \ExtProf \to \calD(\calA)\]
that preserves finite coproducts, and a natural
transformation,
\[\calA[S] \to \calM[S],\]
of condensed anima.

An analytic animated ring
is a pre-analytic animated ring
$(\calA, \calM)$ with the property that
for any object $C \in \calD_{\geq 0}(\calA)$
that is a sifted colimits of objects of the form
$\calM[T]$ for varying extremally disconnected
$T$, the natural map,
\beqn
\label{eqn: analytic animated ring}
\uHom_{\calA}(\calM[S], C) \to \uHom_{\calA}(\calA[S], C),
\eeqn
is an equivalence of condensed anima for all extremally disconnected
profinite sets $S$.
\edefn

To any analytic animated ring $(\calA, \calM)$
we may associate the subcategory of
$\calD_{\geq 0}(\calA)$ of \emph{all} such 
objects $C$ satisfying the condition
\labelcref{eqn: analytic animated ring}
for all profinite sets $S$. 
More precisely, we
have the following definition found in \cite{AnalyticCondensed},
whose notation we have modified for our purposes.

\bdef
\label{def: analytic modules}
Let $\calD_{\geq 0}(\calA, \calM) \subset \calD_{\geq 0}(\calA)$
denote the full $\infty$-subcategory
spanned by all objects
$\calC \in \calD_{\geq 0}(\calA)$ such that
the natural map \labelcref{eqn: analytic animated ring} is
an equivalence of condensed anima for all
extremally disconnected profinite $S$.
\edefn

One may similarly define an abelian subcategory
of $\Mod_{\calA}(\Cond(\Ab))$ for
a given analytic ring.

\bdef
\label{def: discrete analytic modules}
Let $(\calA, \calM)$ be an analytic ring
as defined in \cref{def: analytic ring}.
We denote by $\heart(\calA, \calM) \subset \Mod_{\calA}(\Cond(\Ab))$
the full subcategory of all objects 
$C \in \Mod_{\calA}(\Cond(\Ab))$ such that the map,
\[\uHom_{\calA}(\calM[S], C) \to \uHom_{\calA}(\calA[S], C)\]
is an isomorphism for all
extremally disconnected profinite $S$.
\edefn

The following example of an analytic ring
will be the only one we use in this paper.

\begin{exmp}
\label{exmp: liquid}
Fix $0<p \leq 1$, and consider the
pair $(\R, \calM_{<p})$\footnotemark. This is an
analytic ring by \cite[Theorem 6.5]{AnalyticCondensed},
and the abelian category $\Liq_{p}(\R) := \heart(\R, \calM_{<p})$
is called the category of \emph{$p$-liquid} $\R$-vector spaces.
Likewise, the pair $(\C, \calM_{<p})$ is also an analytic ring,
whose category $\heart(\C, \calM_{<p})$ we also denote $\Liq_p(\C)$.
\end{exmp}

\footnotetext{See Definition 6.3 and page 35 of 
\cite[Lecture VI]{AnalyticCondensed}
for the definition of $\calM_{< p}$.}

\begin{warn}
Note that $\bbR_{\disc}$ is \emph{not}
$p$-liquid as a condensed vector space.
\end{warn}

The following lemma shows that $(\C, \calM_{<p})$
is also an analytic animated ring, under the
inclusion of condensed rings into
condensed animated rings (see
\cref{lem: Top embeds into Cond(Ani)}).

\blem
Suppose that $(\calA, \calM)$ is an analytic
ring. Then $(\calA, \calM)$ is also an analytic
animated ring.
\elem

As such, to each analytic ring $(\calA, \calM)$ we
may associate two categories:
\begin{itemize}
\item[$-$] the abelian category $\heart(\calA, \calM)$ associated to
$(\calA, \calM)$ as an analytic ring,

\item[$-$] and the prestable $\infty$-category 
$\calD_{\geq 0}(\calA, \calM)$ associated to $(\calA, \calM)$
as an analytic \emph{animated} ring.
\end{itemize}

We will use following proposition from \cite{AnalyticCondensed}
to relate these two categories to each other.

\bprop[{\cite[Proposition 12.4]{AnalyticCondensed}}]
\label{prop: analytic ring}
Let $(\calA, \calM)$ be an analytic animated ring. The $\infty$-category
$\calD_{\geq 0}(\calA, \calM)$ is generated under sifted colimits 
by the objects $\calM[S]$ for varying extremally disconnected 
profinite sets $S$, which are compact projective objects 
of $\calD_{\geq 0}(\calA, \calM)$. The full $\infty$-subcategory
	\[\calD_{\geq 0}(\calA, \calM) \subset \calD_{\geq 0}(\calA)\]
is stable under all limits and colimits and admits a left adjoint
\[- \otimes_{\calA} (\calA, \calM): \calD_{\geq 0}(\calA) \to \calD_{\geq 0}(\calA, \calM)\]
characterized as the unique functor
commuting with colimits that
sends $\calA[S]$ to $\calM[S]$.

The $\infty$-category $\calD_{\geq 0}(\calA, \calM)$ is a
prestable. Its heart is 
the full abelian subcategory of $\Mod_{\pi_0 \calA}(\Cond(\Ab))$ generated 
under colimits by $\pi_0 \calM[S]$ for varying $S$. An
object $\calC \in \calD_{\geq 0}(\calA)$ lies in $\calD_{\geq 0}(\calA,\calM)$ if 
and only if all $H_i(\calC)$ lie in $\calD_{\geq 0}(\calA, \calM)^{\heart}$.

If $\calA$ has the structure of an animated condensed 
commutative ring so that $\calD_{\geq 0}(\calA)$ is naturally
a symmetric monoidal $\infty$-category, there is a 
unique symmetric monoidal structure on $\calD_{\geq 0}(\calA, \calM)$
making $- \otimes_{\calA} (\calA, \calM)$ symmetric monoidal.
\eprop

\bcor
\label{cor: vect 1}
Let $(\calA, \calM)$ be an analytic ring.
The heart of the prestable $\infty$-category 
$\calD_{\geq 0}(\calA, \calM)$
is the abelian category of, $\heart(\calA, \calM)$.
Moreover, $\calD_{\geq 0}(\calA, \calM)$
is the connective part of a t-structure,
compatible with filtered colimits,
on a stable presentable
$\infty$-category which we denote
by $\calD(\calA, \calM)$.
\ecor

\bproof
By \cref{prop: analytic ring}, the heart 
of $\calD_{\geq 0}(\calA, \calM)$ is the full subcategory
of condensed $\calA$-modules generated 
under colimits by $\calM[S]$ for varying extremally
disconnected profinite $S$.
But by \cite[Proposition 7.5]{Condensed}, the collection
$\calM[S]$ for $S$ extremally disconnected are a
family of compact projective generators for $\heart(\calA, \calM)$. 
Thus, $\calD_{\geq 0}(\calA, \calM)^{\heart}$ and $\heart(\calA, \calM)$
are two abelian categories generated under
colimits by the same full subcategory
$\scrC^{\cp} \subset \Mod_{\calA}(\Cond(\Ab))$ on objects,
$\{\calM[S]\}_{S \in \ExtProf}$. Moreover,
objects of $\scrC^{\cp}$ are compact projectives in both $\heart(\calA, \calM)$ and
$\calD_{\geq 0}(\calA, \calM)^{\heart}$.
Thus, $\heart(\calA, \calM) \simeq \calD_{\geq 0}(\calA, \calM)^{\heart}$.
The remainder of the claim follows from the
fact that $\calD_{\geq 0}(\calA, \calM)$ is
prestable.
\eproof

The following corollary is immediate
from \cref{prop: analytic ring} and \cref{cor: vect 1}.

\bcor
\label{cor: vect 2}
Let $(\calA, \calM)$ be an analytic ring.
An object $C \in \calD(\calA)$ lies in
$\calD(\calA, \calM)$ if and only
if $H_i(C) \in \heart(\calA, \calM)$ for all $i$.
\ecor

\brem
Clearly, the results of this section
involving $\calD(\calA, \calM)$ are equally
valid for its bounded variations:
$\calD^b(\calA, \calM)$, $\calD^-(\calA, \calM)$,
and $\calD^-(\calA, \calM)$.
\erem

\subsection{The derived category of $\heart(\calA, \calM)$}

\subsubsection{} We fix an analytic ring $(\calA, \calM)$
for the remainder of this section. We let $\Cond(\calA) 
:= \Mod_{\calA}(\Cond(\Ab))$.

\subsubsection{}

\bprop
The category $\heart(\calA, \calM)$ is a Grothendieck
abelian category.
\eprop

\bproof
Recall that a Grothendieck abelian category
is an abelian category that possesses 
arbitrary coproducts, in which filtered
colimits are exact, and which has a single generator\footnotemark.
Observe that $\Cond(\calA)$
is a complete and cocomplete abelian
category with exact filtered colimits by 
\cite[Lemma 18.14.2]{Stacks}
because it is the category of modules over
a sheaf of rings on the topos, $\Cond(\Ab)$.

\footnotetext{More succinctly: an AB5 category with
a generator.}

By \cite[Proposition 7.5]{Condensed}, $\heart(\calA, \calM)$ is a full 
subcategory of $\Cond(\calA)$ that is stable under arbitrary
limits and colimits, meaning that limit and colimits over arbitrary
diagrams in $\Cond(\Ab)$ whose terms lie in $\heart(\calA, \calM)$ also
lie in $\heart(\calA, \calM)$. As mentioned above, $\Cond(\calA)$
admits all colimits and limits, so $\heart(\calA, \calM)$ possesses all colimits,
arbitrary coproduct in particular, and limits, as well.
Since filtered colimits in $\Cond(\calA)$ are exact, 
it suffices to show that the
the natural inclusion $\heart(\calA, \calM) \hook \Cond(\calA)$ is an
exact functor. But this is clear as the inclusion functor preserves
arbitrary limits and colimits\footnotemark.

\footnotetext{Indeed, since $\heart(\calA, \calM) \subset \Cond(\Ab)$ is a full
subcategory, any limit cone in $\heart(\calA, \calM)$ is also
a limit cone in $\Cond(\Ab)$; and since $\heart(\calA, \calM)$
is stable under limits, the limit formed in $\Cond(\Ab)$
is also the limit formed in $\heart(\calA, \calM)$. Ditto for colimits and co-cones.}

It remains to show that $\heart(\calA, \calM)$ has a generator.
By \cite[Theorem 6.5]{AnalyticCondensed}, it is 
generated by compact projective objects,
$\{\calM[S]\}_{S \in \Prof}$. We claim that $\bigoplus_{S \in \Prof}
\calM[S]$ is such a generator. To show that it generates
$\heart(\calA, \calM)$, we must show that for any two distinct morphisms
$f,f': X \to Y$, there exists a morphism $g: \bigoplus_{S \in \Prof} \calM[S]
\to X$, such that $f \circ g \neq f' \circ g$. Since $\{\calM[S]\}_{S \in \Prof}$
are a collection of generators for $\heart(\calA, \calM)$, there exists such
a morphism $g_0: \calM[S_0] \to X$ for some $S_0 \in \Prof$.
Now let $g$ be the morphism induced by $g_0$ and the zero morphism
$0: S \to X$ for all $S \neq S_0$ via the universal property of the
direct sum. It is clear that this satisfies $f \circ g \neq f' \circ g$,
so we are done.
\eproof

As a Grothendieck abelian category, the unbounded
derived $\infty$-category of $\heart(\calA, \calM)$ is a
presentable stable $\infty$-category with a well-behaved
t-structure\footnotemark.

\footnotetext{See \cite[Proposition 1.3.5.9 and 
Proposition 1.3.5.21]{HA}.}

\bprop
\label{prop: D Liq}
There is an equivalence of categories,
\[\calD(\heart(\calA, \calM)) \simeq \calD(\calA, \calM).\]
\eprop

\bproof
It suffices to furnish an equivalence of
Grothendieck prestable $\infty$-categories,
$\calD(\heart(\calA, \calM))_{\geq 0} \simeq 
\calD_{\geq 0}(\calA, \calM)$.

Let $\calD_{\geq 0}(\calA, \calM)^{\cp}$ denote
the full subcategory of compact projective objects
of $\calD_{\geq 0}(\calA, \calM)$. Since 
$\calD_{\geq 0}(\calA, \calM)$ is projectively
generated and admits all small colimits, it is equivalent
to the non-abelian derived category of a minimal model\footnotemark
\footnotetext{See \cite[Definition 2.3.3.1]{HTT} for the
definition of a minimal $\infty$-category. A minimal
model of an $\infty$-category $\scrC$ is a subcategory
of $\scrC$ which is both minimal and equivalent to $\scrC$.} 
of $\calD_{\geq 0}(\calA, \calM)^{\cp}$, which
is closed under finite coproducts by the definition
of analytic animated ring.
The latter category is a $1$-category, which are
always trivially minimal. Thus, there is an equivalence,
\[\calP_{\Sigma}(\calD_{\geq 0}(\calA, \calM)^{\cp}) 
	\simeq \calD_{\geq 0}(\calA, \calM).\]

We remark that $\calD_{\geq 0}(\calA, \calM)^{\cp}$
is an additive $\infty$-category because it contains
the zero object in condensed animated
$\calA$-modules, $\calM[\emptyset] \simeq 0$,
and is closed under finite finite biproducts; it is also
small because it is a full subcategory of a locally
small category ($\calD_{\geq 0}(\calA, \calM)$)
whose collection of objects is indexed by 
the small category $\ExtProf$.
As such, $\calD_{\geq 0}(\calA, \calM)$
is the non-abelian derived category of a small
additive $\infty$-category, so by \cite[Remark C.1.5.9]{SAG}
and \cite[Proposition C.5.3.4]{SAG},
it is a $0$-complicial complete 
Grothendieck prestable $\infty$-category.
This combination of adjectives implies,
by \cite[Corollary C.5.9.7]{SAG}, that
the inclusion of the heart,
$\heart(\calA, \calM) \hook \calD_{\geq 0}(\calA, \calM)$,
extends to an equivalence,
\beqn
\label{eqn: complete derived}
\wih{\calD}(\heart(\calA, \calM))_{\geq 0} \xrightarrow{\simeq} \calD_{\geq 0}(\calA, \calM),
\eeqn
of $\infty$-categories.

By \cite[Proposition C.5.9.2]{SAG}, 
the completion functor on $\Groth_{\infty}^{\lex}$ 
restricts to an equivalence of categories,
\[\Groth_{\infty}^{\ch, \lex} 
	\xrightarrow[\simeq]{\wih{(-)}} \Groth_{\infty}^{\comp, \lex},\]
between anticomplete Grothendieck prestable
$\infty$-categories and complete Grothendieck
prestable $\infty$-categories.
By the universal property of the
unseparated derived prestable $\infty$-category
(\cite[Corollary C.5.8.9]{SAG}), there
exists a map
\beqn
\label{eqn: unseparated derived}
\widecheck{\calD}(\heart(\calA, \calM))_{\geq 0}
	\to \calD_{\geq 0}(\calA, \calM)
\eeqn
extending the inclusion of the heart.
By unraveling the definitions, 
we see that this map (which is
also the map in \cite[Corollary C.5.8.11]{SAG}) is sent under
$\wih{(-)}$ to precisely the functor
\labelcref{eqn: complete derived}. Since
$\wih{(-)}$ is an equivalence of categories,
it follows that \labelcref{eqn: unseparated derived}
is also an equivalence.

Thus, we have that $\widecheck{\calD}(\heart(\calA, \calM))_{\geq 0} \simeq 
\wih{\calD}(\heart(\calA, \calM))_{\geq 0}$ by \labelcref{eqn:
complete derived} and \labelcref{eqn: unseparated derived}.
We claim that this implies that $\widecheck{\calD}(\heart(\calA, \calM))_{\geq 0}
\simeq \calD(\heart(\calA, \calM))$. 
By \cite[Theorem C.5.4.9]{SAG},
$\calD(-)_{\geq 0}: \Groth_{\ab}^{\lex} \to \Groth_{\infty}^{\lex, \sep}$ 
is the left adjoint
to the functor of restriction
to the heart. That is,
\[\LFun^{\lex}(\calD(\scrA)_{\geq 0}, \scrC) \simeq \LFun^{\lex}(\scrA, \scrC^{\heart}),\]
where $\LFun^{\lex}(-,-)$ denotes the morphisms
of $\Groth_{\infty}^{\lex, \sep}$ and $\Groth_{\ab}^{\lex}$\footnotemark.
\footnotetext{This is Lurie's notation in \cite[Appendix C]{SAG}.
When both arguments of $\LFun(-,-)$ are presentable categories (such as is
the case for Grothendieck prestable categories), $\LFun$ denotes
functors that preserve small colimits. The superscript ``$\lex$"
denotes those functors which further preserve finite limits.}
On the other hand, $\widecheck{\calD}(-)_{\geq 0}$ is
the left adjoint in $\Groth_{\infty}^{\lex}$ to the
functor of restriction to the heart by \cite[Corollary C.5.8.9]{SAG}. 
The functor $\scrC \mapsto \scrC^{\sep}$ taking a prestable
category to its separable quotient is left adjoint to
the inclusion, $\Groth_{\infty}^{\lex, \sep} \to \Groth_{\infty}^{\lex}$,
by \cite[Corollary C.3.6.2]{SAG}. By the essential uniqueness
of adjoint functors, we therefore have 
\[(\widecheck{\calD}(-)_{\geq 0})^{\sep}
	\xrightarrow{\simeq} \calD(-)_{\geq 0}.\]
But now observe that
$\widecheck{\calD}(-)_{\geq 0}$ is complete,
and therefore already separated, so
$\widecheck{\calD}(-)_{\geq 0} \simeq (\widecheck{\calD}(-)_{\geq 0})^{\sep}$.

All together, we obtain:
\begin{align*}
\calD(\heart(\calA, \calM))_{\geq 0} 	&\simeq \widecheck{\calD}(\heart(\calA, \calM))_{\geq 0} \\
							&\simeq \wih{\calD}(\heart(\calA, \calM))_{\geq 0} \\
							&\simeq \calD_{\geq 0}(\calA, \calM),
\end{align*}
which completes the proof.
\eproof

\bprop
\label{prop: internal hom}
Both $\heart(\calA, \calM)$ and
$\calD(\calA, \calM)$ are closed symmetric
monoidal categories with respect to the
symmetric monoidal structures specified in
\cite[Proposition 7.5]{Condensed}
and \cite[Proposition 12.4]{AnalyticCondensed}, respectively.
Moreover, $\otimes_{\calD^-(\calA, \calM)}$
is the left derived functor of
$\otimes_{\heart(\calA, \calM)}$, and $\uHom_{\calD^+(\calA, \calM)}$
is the right derived functor of $\uHom_{\heart(\calA, \calM)}$.
\eprop

\bproof
By \cref{prop: analytic ring}, the left adjoint,
$- \otimes_{\calA} (\calA, \calM)$ is symmetric
monoidal. Thus, the tensor product in
$\calD(\calA)$ has a right adjoint,
denoted $\uHom_{\calD(\calA)}(-,-)$, the tensor
product in $\calD(\calA, \calM)$ also has a right adjoint,
given by $\uHom_{\calD(\calA)}(-,-) \otimes_{\calA} 
(\calA, \calM)$.
The same argument shows
that the application of the left adjoint
to the inclusion $\heart(\calA, \calM)
\hook \Mod_{\calA}(\Cond(\Ab))$
(the liquidification if $(\calA, \calM)$
were $(\C, \calM_{<p})$)
to $\uHom_{\Cond(\calA)}(-,-)$
gives the internal hom in $\heart(\calA, \calM)$.

The claim that $\otimes_{\calD^-(\calA, \calM)}$
is the left derived functor of $\otimes_{\heart(\calA, \calM)}$
follows immediately from \cref{prop: D Liq} and the
uniqueness of left derived functors, 
Indeed, if $P \in \heart(\calA, \calM)$, then
$P \otimes -$ determines a right t-exact functor\footnotemark
\footnotetext{Because it tautologically restricts
to a symmetric monoidal functor on 
$\calD_{\geq 0}(\calA, \calM)$.},
\[P \otimes -: \calD^-(\heart(\calA, \calM)) \to \calD^-(\calA, \calM),\]
that restricts to $P \otimes_{\heart(\calA, \calM)} -$ on the hearts.
This functor is right exact as it commutes with colimits;
so, by \cite[Theorem 1.3.3.2]{HA}, it is
the left derived functor of $\otimes_{\heart(\calA, \calM)}$
(noting that $\calD^-(\calA, \calM)$ is left-complete).

The claim that $\uHom_{\calD^+(\calA, \calM)}$
is the right derived functor of $\uHom_{\heart(\calA, \calM)}$ now
follows immediately.
\eproof

\section{Functional calculus via condensed mathematics}

\subsubsection{} Fix $0<p \leq 1$. 

\subsubsection{}
We now specialize the results of the
previous section to the analytic ring
$(\C, \calM_{<p})$ recalled
in \cref{exmp: liquid}.

The theory of $p$-liquid vector spaces,
i.e.\ objects of $\Liq_p(\C)$, is an excellent
framework in which to do functional analysis.
Many of the most commonly encountered
types of topological
vector spaces are $p$-liquid
vector spaces,
which, as seen above, enjoy
great homological properties as
a category.
Banach spaces, viewed as $p$-liquid
vector spaces, will be our main
objects of interest.

\subsection{The category of Banach spaces}
We consider the following two categories
of topological vector spaces:
\begin{enumerate}[(i)]

\item Let $\TVect^{\comp}_{\lc}$ denote the category
of complete\footnotemark locally convex topological $\C$-vector spaces
with morphisms given by continuous linear operators.
\footnotetext{i.e.\ the underlying topological
space is completely metrizable.}

\item Let $\Ban$ denote the category of
complex Banach spaces with morphisms given
by bounded linear maps

\end{enumerate}

There is an obvious forgetful functor 
$\oblv: \Ban \to \TVect^{\comp}_{\lc}$
sending a Banach space to its underlying
locally convex topological vector space.
The functor $\oblv$ is fully faithful;
every continuous linear operator between
Banach spaces is bounded.

It is shown in \cite{AnalyticCondensed} 
that $\Ban$ embeds fully faithfully into $p$-liquid
vector spaces as follows.
Complete locally convex topological
spaces are compactly generated, so they embed
into condensed sets (see discussion after 
\cite[Definition 3.1]{AnalyticCondensed}). Moreover
they are $\calM$-complete in the sense of 
\cite[Definition 4.1]{AnalyticCondensed} by
\cite[Proposition 3.4]{AnalyticCondensed},
so they are $p$-liquid by the 
discussion at the start of \cite[\S VI]{AnalyticCondensed}.
The embedding of Banach spaces into $p$-liquid
spaces factors as
\beqn
\label{eqn: Ban embedding}
\Ban \hook \TVect^{\comp}_{\lc} \hook \Liq_p(\C).
\eeqn

Using this embedding, we make the following
definition.

\bdef
Let $\calD_{\Ban}(\C, \calM_{<p})$ denote
the full subcategory of $\calD(\C, \calM_{<p})$
spanned by objects whose homology spaces lie in
the essential image of the embedding $\Ban
\hook \Liq_p(\C)$.
\edefn

\subsection{Classical functional calculus}
We briefly review the classical theory of
functional calculus on Banach spaces.
Our reference for the content of this section
is the introductory paragraphs of \cite[Chapter 10]{RudinFun}.

\subsubsection{Motivation of functional calculus}
Given a Banach space $V \in \Ban$, and a bounded
linear operator $T: V \to V$, the symbol ``$T^n$"
has a clear, unambiguous meaning for each non-negative
integer $n$. For $n>0$, $T^n$ is the $n$-fold composition
of $T$ with itself, and $T^0 = I$.
Likewise, given a polynomial $f(z)$ with complex coefficients,
the symbol ``$f(T)$" has a clear meaning as well. The
purpose of functional calculus is to try and extend the
definition of symbols ``$f(T)$" to include functions $f$ which
are holomorphic on an open neighborhood of $\C$ containing
the spectrum of $T$.

\subsubsection{}
This is easily achieved in the classical setting
by observing that the space of bounded endomorphisms,
$\End_{\Ban}(V)$, is itself a Banach space under the
strong operator norm. In fact, it is a Banach \emph{algebra} under
the composition of bounded operators.

The following is an easy result 
in classical functional analysis that
uses the fact that polynomials are dense
inside the space the space of holomorphic
functions.

\bprop
\label{prop: classical functional calculus}
Let $A$ be a Banach algebra,
$a \in A$ be an element, and $f(z)$ 
be an entire function of one
complex variable. Then there is
a unique element $f(a)$ extending
the definition of $f(a)$ from polynomials
to entire functions in a continuous way.
\eprop

In the sequel, we will often refer to
the endomorphism $f(T) \in \End_{\Ban}(V)$
obtained ``from classical functional calculus."
By this, we mean the endomorphism obtained
from the application of $f$ to $T$ as an element
of the Banach algebra $\End_{\Ban}(V)$ in the
sense of \cref{prop: classical functional calculus}.

\subsection{Liquid functional calculus}
\label{ssec: liquid functional calculus}

Suppose given $E \in \calD(\C, \calM_{<p})$
such that the internal hom, 
$\End(E)$, lies in the subcategory $\calD_{\Ban}(\C, \calM_{<p})$.
The goal of this section is to make sense of
$f(T)$ for any $T \in \End(E)$ and 
any entire function $f$
in such a way that the maps
$H_i(-)$ induced by $f(T)$ 
coincide with the operators
obtained from classical functional calculus.

\subsubsection{}

Let $f(z) = \sum_{i=0}^{\infty} a_i z^i$
be an entire function of a single variable.
Since $\End(E) \in \calD(\C, \calM_{<p})$,
the following morphism of condensed 
$\un{\C}$-module mapping spectra\footnotemark
\footnotetext{See \cref{lem: internal hom} for proof of the 
existence of internal hom in $\calD(\C)$ and 
the discussion above for remarks about the
underlying condensed spectra of an object in
$\calD(\C)$.} 
is equivalent,
\[\uHom_{\calD(\C)}(\calM_{< p}(S), \End(E)) \to \uHom_{\calD(\C)}(\C[S], \End(E)),\]
for any profinite set $S$.
The latter, in turn, is equivalent as a condensed
spectrum to $\uHom_{\Cond(\Sp)}(S, \End(E))$, because
$\C[S]$ is the free object on $S$ in $\calD(\C)$ (i.e.\ is the left adjoint to the
forgetful functor $\calD(\C) \to \Cond(\Sp)$).

\bdef
\label{def: f as a measure}
Observe that $S := \bbN \cup \{\infty\}$ is
a profinite set. Let $f$ be as above. 
The element $\mu_f \in \C[S]$ is
defined by the assignment
\begin{align*}
n &\mapsto \sqrt{a_n} \\
\infty &\mapsto 0. 
\end{align*}
\edefn

\bprop
\label{prop: mu is a measure}
The measure $\mu_f$ is an element
of $\calM_{<p}(S)$.
\eprop

\bproof

Consider $S$ as the projective limit
of the (totally ordered) finite sets $S_i := \{0, \ldots, i, \infty\}$,
where the transition maps
$f_{j>i}: S_j \to S_i$ for the
directed diagram are given by
\[
f_{j>i}(x)=
    \begin{cases}
        x & \text{if } 0 \leq x \leq i\\
        \infty & \text{if } x > i
    \end{cases}.
\]
Recall that for any finite set $T$,
\[\R[T]_{\ell^p \leq c} := \{(a_t)_{t \in T} \in \R[T] | \sum_t |a_t|^p \leq c\},\]
where $c > 0$, and that for a profinite
set, $T = \varprojlim_i T_i$,
\[\calM_p(T) := \bigcup_{c>0} \varprojlim_i \R[T_i]_{\ell^p \leq c}.\]

We claim that $\mu_f \in \calM_p(S)$ for any $0 < p \leq 1$.
Indeed, observe that $\R[S] \simeq \varprojlim_i \R[S_i]$
and consider the element $\mu_i \in \R[S_i]$
given by the assignment
$n \mapsto \sqrt{a_n}$, $\infty \mapsto 0$.
Clearly, $\mu_f = (\ldots, \mu_i, \ldots, \mu_0)
\in \varprojlim_i \R[S_i]$ under the identification
of $\R[S]$ with the projective limit $\varprojlim_i \R[S_i]$.

It now suffices to show that $\mu_i \in \R[S_i]_{\ell^p \leq c_0}$,
for some $c_0 >0$ independent of $i$.
Since $f(z)$ is holomorphic, its power
series representation $\sum_{i=0}^{\infty} a_i z^i$
is absolutely convergent for any value of $z$.
In particular, by \cref{lem: holomorphic ratio test}
and \cref{lem: ratio test for powers} stated
at the end of this section, we obtain
that the series $\sum_{i=0}^{\infty} |a_i|^p$ is
convergent. Applying \cref{lem: ratio test for powers}
a second time to this latter series, 
we see that $\sum_{i=0}^{\infty} \left|\sqrt{a_i} \right|^p$
is convergent. Since $\sum_{i=0}^{m} \left|\sqrt{a_i} \right|^p
\leq \sum_{i=0}^{\infty} \left|\sqrt{a_i} \right|^p$ for any choice
of $m$, we may take $c_0 := \sum_{i=0}^{\infty} |\sqrt{a_i}|^p$. Thus, 
$\mu_f \in \calM_p(S)$.

Finally, recall that $\calM_{<p}(T) := \bigcup_{p'<p} \calM_{p'}(T)$
for any profinite set $T$. Since $\mu_f \in \calM_p$(S) for any
fixed $0 < p' \leq 1$, it belongs in particular to $\calM_{p'}(S)$
for some $p' <p$, and therefore belongs to $\calM_{<p}(S)$.
\eproof

\Cref{prop: mu is a measure} allows us
to make sense of the application of $f$ to
an endomorphism of $E$.
Let $T \in \End(E)$ be an endomorphism
of $E$, and denote by $\varphi_{T,f}: S \to \End(E)$
the map of condensed spectra given by
the assignment, 
\begin{align*}
n &\mapsto \sqrt{a_n} \cdot T^n \\
\infty &\mapsto 0,
\end{align*}
where we note that this
indeed defines a map of condensed
spectra by the fully faithful
embedding of spectra into
condensed spectra.
Since $\End(E)$ is $p$-liquid for any $0<p \leq 1$,
\[\uHom_{\C}(\calM_{<p}(S), \End(E)) \xrightarrow{\simeq} \uHom_{\Cond(\Sp)}(S, \End(E))\]
is a equivalence of condensed spectra.
Thus, $\varphi_{T,f}$ extends to a map
$\wit{\varphi_{T,f}}: \calM_{<p}(S) \to \End(E)$.

\bdef
\label{def: liquid functional calculus}
Define $f(T) \in \End(E)$ to be the
endomorphism,
\[f(T) : = \wit{\varphi_{T,f}}(\mu_f),\]
given by evaluating
$\wit{\varphi_{T,f}}$ on the measure
$\mu_f \in \calM_{<p}(S)$.
\edefn

\brem
\label{rem: classical = spectral p-liquid}
Observe that if $E \in \Liq_p(\C)
\subset \calD(\C, \calM_{<p})$, the spaces of morphisms,
$\Hom_{\calD(\C)}(\calM_{<p}(S), E)$ and 
$\Hom_{\Cond(\Sp)}(S,E)$ are discrete. Thus,
if $E$ is $p$-liquid as a $\C$-module spectrum,
it is $p$-liquid as a $\C$-vector space in the
sense of \cite[Theorem 6.5]{AnalyticCondensed}.
Moreover, the extension $\wit{f}: \calM_{<p}(S) \to E$
of a map $f: S \to E$ is the unique extension
guaranteed by \cite[Theorem 6.5]{AnalyticCondensed},
as well.
\erem

\subsection{Comparison of liquid to classical functional calculus}
It remains to show that, if $E$ is a Banach
space viewed a condensed $\C$-vector space, 
the notation $f(T)$
is compatible with its meaning coming
from classical functional analysis.

Given a profinite set $S$, let $\calM(S)$
be the real vector space of signed Radon
measures on $S$\footnotemark. 
\footnotetext{See \cite[Exercise 3.3]{Condensed}
for more context.}
Let $W$ be an arbitrary Banach
space. By \cite[Proposition 3.4]{Condensed},
$W$ is $\calM$-complete in the sense that
any map $f: S \to W$ from a profinite set $S$
extends uniquely to a map of topological vector
spaces,
\begin{align*}
\overline{f}: \calM(S) &\to W \\
\mu &\mapsto \int_S f d\mu
\end{align*}
where the right-hand side
is the integral of a Banach space-valued
function on a compact measure space
(see \cite[Theorem 3.27]{RudinFun}).

Every $\calM$-complete vector space
is $p$-liquid for any fixed $0<p \leq 1$; indeed,
$\calM_{<p}(S) \subset \calM(S)$
for any such $p$. Moreover, since the extension
of $f$ to a function $\wit{f}: \calM_{<p}(S) \to W$
is unique by the $p$-liquid property of $W$,
it follows that $\wit{f} = \over{f}|_{\calM_{<p}(S)}$.

With the preceding discussion in mind, we state
and prove the following lemma, which claims
that a convergent series (equivalently sequence) 
of terms in $W$ can be viewed as the evaluation
of $\wit{f}$ on a particular measure.

\blem
\label{lem: extension given by integration}
Suppose that $W$ is a complex Banach space, and let 
$\sum_{i=0}^{\infty} a_i x_i$, $x_i \in W$, be
a convergent series in $W$ such that 
$\sum_{i=0}^{\infty} a_i$
passes the ratio test.

Let $f: S := \bbN \cup \{\infty\} \to W$ be the
continuous function given by the assignment
$n \mapsto \sqrt{a_n} \cdot x_n$, $\infty \mapsto 0$,
and let $\mu$ be the signed Radon measure on
$S$ given by $\mu(n) = \sqrt{a_n}$. Then,
\beqn
\label{eqn: series in Banach space}
\sum_{i=0}^{\infty} a_i x_i = \wit{f}(\mu).
\eeqn


\elem

\bproof
Firstly, we note that, by \cref{lem:
ratio test for powers}, $\mu \in \calM_{<p}(S)$
for any fixed $p$, so the right-hand side
of \labelcref{eqn: series in Banach space}
is well-defined.

By the preceding discussion,
we have that
$\wit{f} = \over{f}|_{\calM_{<p}(S)}$, so
\begin{align*}
\wit{f}(\mu) &= \over{f}(\mu) \\ 
			&= \int_S f d\mu \\
			&= \sum_{i=0}^{\infty} a_i x_i,
\end{align*}
by the definition of the integral, $\int_S f d\mu$.
\eproof

In particular, suppose $W := \End(V)$
is the Banach algebra of bounded linear operators
on a given Banach space $V$.
As a Banach space, $\End(V)$ is $p$-liquid, so the
map $\varphi_{T,f}: S \to \End(V)$ given by $n \mapsto 
\sqrt{a_n} \cdot T^n$, $\infty \mapsto 0$
for $T \in \End(V)$ extends uniquely to a map,
\[\wit{\varphi_{T,f}}: \calM_{<p}(S) \to \End(V).\]
The following corollary follows immediately
from \cref{lem: extension given by integration}
and \cref{rem: classical = spectral p-liquid}.

\bcor
\label{cor: liquid-classical comparison}
Let $\mu_f \in \calM_{<p}(S)$ be as
in \cref{prop: mu is a measure}. Then
\[\wit{\varphi_{T,f}}(\mu_f) = f(T),\]
where the right-hand side denotes the
operator given by the convergent sum
$\sum_{i=0}^{\infty} a_i T^i$, in the strong
operator norm.
\ecor

\subsection{Passing to homology}

Let $E \in \calD(\C, \calM_{<p})$
and $T \in \End(E)$ be as before.
Then $T$ induces endomorphisms of the
homology groups of $E$ (with respect
to the t-structure on $\calD(\C, \calM_{<p})$), 
via the homology functor,
$H^i: \calD(\C, \calM_{<p}) \to \Liq_p(\C)$.
Namely, for each $i \in \bbZ$,
we obtain a map of condensed
liquid vector spaces
\[H_i(-): \End(E) \to \End(H_i(E))\]
by functoriality.
Since it was assumed $H^i(E)$
lies in the image of Banach spaces
inside of $p$-liquid vector spaces,
these endomorphisms are bounded
operators on complex Banach spaces,
so $f(H_i(T))$ for an entire 
function $f$ is easily defined classically.
The following proposition asserts that
the condensed functional calculus we
have defined recovers the classical
functional calculus upon taking homology.

\bprop
\label{prop: liquid-classical comparison}
The induced map on homology,
\[H_i(f(T)): H_i(E) \to H_i(E),\]
are given by the Banach space
endomorphism obtained via 
classical functional calculus
by applying $f$ to
\[f(H_i(T)): H_i(E) \to H_i(E).\]
\eprop

\bproof
Consider the map $\varphi_{T,f,i}: S := \bbN \cup \{\infty\} \to \End(H_i(E))$
given by the assignment, 
\begin{align*}
n 	&\mapsto \sqrt{a_n} \cdot H_i(T)^n \\
\infty &\mapsto 0.
\end{align*}
The target of this map is a $p$-liquid vector space, 
so $\varphi_{T,f,i}$ extends uniquely
to a map $\wit{\varphi_{T,f,i}}: \calM_{<p}(S) \to \End(H_i(E))$.
By \cref{cor: liquid-classical comparison}, $\wit{\varphi_{T,f,i}}(\mu_f)$
is the endomorphism $f(H_i(T))$ given by classical functional
calculus.

Thus, it remains to show that
\[H_i(\wit{\varphi_{T,f}}(\mu_f)) = \wit{\varphi_{T,f,i}}(\mu_f).\]
In order to see this, we show that
$H_i \circ \wit{\varphi_{T,f}}$ is an
extension to $\calM_{<p}(S)$ of
the map $\varphi_{T,f,i}$, since such
an extension is unique by the
definition of $p$-liquid vector space.
Pictorially, we have
\[\begin{tikzcd}
	S & {\End(E)} & {\End(H_i(E))} \\
	{\calM_{<p}(S)}
	\arrow["{H_i}", from=1-2, to=1-3]
	\arrow["{\varphi_{T,f}}", from=1-1, to=1-2]
	\arrow[""', hook', from=1-1, to=2-1]
	\arrow["{\wit{\varphi_{T,f}}}"{pos=0.3}, from=2-1, to=1-2]
	\arrow["{H_i(\wit{\varphi_{T,f}})}"', from=2-1, to=1-3]
\end{tikzcd}\]
from which we see that $(H_i \circ \wit{\varphi_{T,f}})|_{S}
= H_i \circ \varphi_{T,f}$.

Since $H_i: \calD(\C, \calM_{<p}) \to 
\calD(\C, \calM_{<p})^{\heart} \simeq \Liq_p(\C)$ 
is an additive functor of categories enriched
over condensed $\C$-vector spaces, $H_i: \End(E)
\to \End(H_i(E))$ is $\C$-linear.
As such, we have the equality,
\[H_i(c \cdot T^n) = c \cdot H_i(T)^n\]
for any natural number
$n \in \bbN$ and complex scalar $c \in \bbC$.
Computing, we obtain
\begin{align*}
(H_i \circ \varphi_{T,f})(n) 	&= H_i(\sqrt{a_n} \cdot T^n) \\
					&= \sqrt{a_n} \cdot H_i(T)^n \\
					&= \varphi_{T,f,i}(n).
\end{align*}
To conclude, we observe that $(H_i \circ \varphi_{T,f})(\infty) = 0$.
\eproof

\subsection{Postponed lemmas}

We state and prove two lemmas on
the ratio test for convergence of series that
were used in the proof of \cref{prop: mu is a measure}.

\blem
\label{lem: holomorphic ratio test}
Let $f(z) = \sum_{i=0}^{\infty} a_i z^i$ be
any entire function of a single complex variable.
Then the series $\sum_{i=0}^{\infty} |a_i|$ passes
the ratio test for convergence.
\elem

\bproof
For fixed $n \in \bbN$ and $z \in \C$, we have the following equality,
\[\left| \frac{a_{n+1}z^{n+1}}{a_n z^n} \right| = \left| \frac{a_{n+1}}{a_n}\right| \left| z\right|.\]
Taking the limits as $n \to \infty$, we obtain
\beqn
\label{eqn: holomorphic ratio test}
\lim_{n \to \infty} \left| \frac{a_{n+1}z^{n+1}}{a_n z^n} \right| = \left( \lim_{n \to \infty} \left| \frac{a_{n+1}}{a_n}\right| \right) |z|.
\eeqn
Since $\sum_{i=0}^{\infty} a_i z^i$ is absolutely
convergent for all $z \in \C$, the left-hand
side of \labelcref{eqn: holomorphic ratio test}
converges to a value $\leq 1$ for any value of $z$, for otherwise
$\sum_{i=0}^{\infty} a_i z^i$ would diverge by the ratio test.
Now, take any $z = z_0$ such that $|z_0|>1$ to obtain that
\[\lim_{n \to \infty} \left| \frac{a_{n+1}}{a_n}\right| <1,\]
which was to be shown.
\eproof

\blem
\label{lem: ratio test for powers}
If $\sum_{i=0}^{\infty} a_i$ passes the ratio
test, then $\sum_{i=0}^{\infty} {a_i}^r$ does
too, for any fixed $0<r$.
\elem

\bproof

We compute,
\[\left| \frac{{a_{n+1}}^r}{{a_n}^r} \right| = \left| \left( \frac{a_{n+1}}{a_n} \right)^r \right| 
								= \left| \frac{a_{n+1}}{a_n} \right|^r.\]
Since $(-)^r$ is a continuous
function of a complex variable,
\beqn
\label{eqn: ratio test for powers}
\lim_{n \to \infty} \left( \left| \frac{a_{n+1}}{a_n} \right|^r \right) =
	 \left( \lim_{n \to \infty} \left| \frac{a_{n+1}}{a_n} \right| \right)^r.
\eeqn
By assumption, $\lim_{n \to \infty} \left| \frac{a_{n+1}}{a_n} \right| < 1$,
so the right-hand side of \labelcref{eqn: ratio test for powers}
is also $<1$, proving the lemma.
\eproof

\section{Functional calculus in $\Perf$}
\label{sec: Functional calculus in Perf}
The results of the previous section
allow us to define things such as
the exponential of an endomorphism
in $\Vect^b(\C)$ of an object with finite
dimensional homology spaces (i.e.\ a
perfect complex), as we show in this
section.

\blem
$\Ban$ is an
additive category.
\elem

\bproof
Clearly, the hom-sets in $\Ban$ are
abelian groups: the point-wise
sum of any two bounded linear operators is again
a bounded linear operator. The composition
of morphism in $\Ban$ is also clearly bilinear,
so $\Ban$ is enriched over $\Ab$.
It remains to check that $\Ban$ admits
finite products. The product of two Banach
spaces $V_1$ and $V_2$ is canonical a Banach
space under several equivalent norms,
such as the max norm: 
$\|\boldsymbol{\cdot}\|_{V_1 \times V_2}
:= \mathrm{max}(\|\boldsymbol{\cdot}\|_{V_1}, 
\|\boldsymbol{\cdot}\|_{V_2})$. 
\eproof

Since $\Ban$ is additive, we can consider
its category of chain complexes, $\Ch(\Ban)$,
which is also additive. 

\brem
Neither $\Ban$ nor
$\Ch(\Ban)$ are abelian.
In fact, one can see the theory of
$p$-liquid vector spaces as a
remedial solution to this fact.
\erem

Notably, this embedding is \emph{additive}. The functor
$\Ban \hook \TVect^{\comp}_{\lc}$ obviously preserves
finite products; we claim that $\TVect^{\comp}_{\lc} \hook \Liq_p(\C)$ does too.
In fact, we prove a stronger claim.

\blem
\label{lem: preserve limits}
The fully faithful embedding $\TVect^{\comp}_{\lc}
\hook \Liq_p(\C)$
preserves limits\footnotemark.
\elem

\footnotetext{cf.\ \cite[Remark 1.8]{Condensed}}

\bproof
By \cite[Proposition 1.7]{Condensed}, the embedding
$\Top_{\kappa} \hook \Cond_{\kappa}$ of $\kappa$-compactly
generated topological spaces into $\kappa$-condensed sets
is a right adjoint functor, so preserves limits. Thus, it suffices
to show that the forgetful functors, $\oblv_1: \TVect^{\comp}_{\lc}
\to \Top_{\kappa}$ and $\oblv_2: \Cond_{\kappa}(\C) \to \Cond_{\kappa}$,
reflect limits.

We begin with $\oblv_1$.
By \cite[Proposition 3.3]{nlab},
a conservative functor reflects any limits that
exist in the domain and which it preserves.
The functor $\oblv_1$ is conservative by
\cite[2.12 Corollaries (a) and (b)]{RudinFun},
so it suffices to show that $\oblv_1$ preserves
all limits.
By \cite[Lemma 4.14.11]{Stacks}, it 
suffices to check that $\oblv_1$
preserves products and equalizers.
Products are preserved under the forgetful functor,
since the underlying topological space of a product
in $\TVect^{\comp}_{\lc}$ has underlying topological
space given by the Cartesian product of the underlying 
topological spaces of each of the factors. 
Since $\TVect^{\comp}_{\lc}$ is an additive
category, it suffices just to check that equalizers
of the form,
\[\begin{tikzcd}
	V & W
	\arrow["f", curve={height=-6pt}, from=1-1, to=1-2]
	\arrow["0"', curve={height=6pt}, from=1-1, to=1-2],
\end{tikzcd}\]
i.e.\ kernels, are preserved. But the kernel
of $f$ in $\TVect^{\comp}_{\lc}$ is given
by the vector subspace $f^{-1}(0)$, since this is
automatically closed, and therefore complete.
We now conclude by noting that the 
underlying topological space of $f^{-1}(0)$
is precisely the fiber product $\ast \,{_0\times_f} V$
taken in $\Top_{\kappa}$.

We employ the same strategy to show
that $\oblv_2$ reflects limits. The forgetful functor
$\oblv_2$ is conservative and preserves
limits for general sheaf-theoretic
reasons: $\Cond_{\kappa}(\C)$ is the category of
modules over a sheaf of rings on a small site ($\Prof_{\kappa}$),
so the forgetful functor to the category of $\Set$-valued
sheaves is conservative by \cite[Corollary 3.4.4.6]{HA},
and preserves limits by \cite[Corollary 3.4.3.2]{HA}.
Thus, $\oblv_2$ reflects limits, which concludes the proof.
\eproof

Since all norms on a finite dimensional
complex vector space are equivalent and complete, and
any linear map between such objects is bounded,
there is an inclusion $\Vect_{\C}^{f.d.} 
\hook \Ban$ where we view $\Ban$ as a full subcategory of
$\TVect^{\comp}_{\lc}$. It is not hard to
see that this inclusion is an additive functor,
so composing with \labelcref{eqn: Ban embedding},
we obtain the fully faithful, additive functor,
\beqn
\label{eqn: Vect embedding}
\un{(-)}:\Vect_{\C}^{f.d.} \hook \Liq_p(\C),
\eeqn
of abelian categories.

\begin{warn}
The embedding \labelcref{eqn: Vect embedding}
is \emph{not} the restriction of the embedding
$\Vect_{\C} \hook \Cond(\Ab)$ given by
viewing $\C$ as a discrete topological ring
and sending an object in $\Vect_{\C}$ to a
$\un{\C_{\disc}}$-module in condensed
abelian groups.
\end{warn}

\blem
\label{lem: exact}
The embedding $(\Vect_{\C}^{f.d.}, \otimes) 
\xhookrightarrow{\un{(-)}} (\Liq_p(\C), \otimes_{\Liq_p})$
is exact and symmetric monoidal.
\elem

\bproof
To show that $\un{(-)}$ is exact, we must 
show that it preserve finite limits and colimits.
We have that $\un{(-)}$ preserves finite limits
by \cref{lem: preserve limits}, so we need only
show that it preserves finite colimits. 
For this, it suffices by \cite[Lemma 4.14.12]{Stacks},
to show that $\un{(-)}$ preserves cokernels
and finite coproducts.

Recall that $\un{(-)}$ sends
an object $V \in \Vect_{\C}^{f.d.}$ to the
sheaf given by the assignment
\[S \mapsto \Hom_{\Top}(S, V),\]
for $S$ an extremally disconnected
profinite set, where $V$ is considered
with its canonical Euclidean topology. 
For objects $V, W \in \Vect_{\C}^{f.d.}$,
the presheaf given by
\[S \mapsto \un{V}(S) \oplus \un{W}(S)\]
is actually a sheaf because sheafification
commutes with colimits as a left adjoint functor. 
This sheaf is the coproduct $\un{V} \oplus \un{W}$.
Then we have
\begin{align*}
\un{V} \oplus \un{W}(S) &= \un{V}(S) \oplus \un{W}(S) \\
						&= \un{V}(S) \times \un{W}(S) \\
						&= \Hom_{\Top}(S, V) \times \Hom_{\Top}(S, W) \\
						&= \Hom_{\Top}(S, V \times W) \\
						&= \un{(V \times W)}(S) \\
						&= \un{(V \oplus W)}(S),
\end{align*}
where $\un{V}(S) \oplus \un{W}(S) = \un{V}(S) \times \un{W}(S)$
because $\oplus$ and $\times$ here are
taken in $\Ab$, where they coincide.
Thus, $\un{(-)}$ preserves coproducts.
The proof that $\un{(-)}$ preserves
cokernels is similar, using that cokernels
are colimits.

Finally, to see that $\un{(-)}$ is symmetric
monoidal, we note that $\Vect_{\C}^{f.d.}$ is
a full subcategory of the category of \emph{dual nuclear Fr\'echet}
spaces, which embed in the obvious way into $\Liq_p(\C)$
as a subset of \emph{nuclear} spaces\footnotemark. 
\footnotetext{See \cite[Lecture VIII]{ComplexCondensed}
for the definitions of dual nuclear Fréchet space and
nuclear space in this context.}
Now we turn to
Scholze's initial post on the Xena Project blog (\cite{liquidblog}) 
putting forth the mathematical formalization
challenge known as the ``Liquid Tensor Experiment,"
in which he mentions that $\otimes_{\Liq_p}$ agrees with
the usual completed tensor product on nuclear spaces,
so we are done.
\eproof

Observe that both $\Vect_{\C}^{f.d.}$ and
$\Liq_p(\C)$ have enough projectives. 
Indeed, all objects in $\Vect_{\C}^{\heart}$ are
projective, and since $\Vect_{\C}^{f.d.}$ is a full
subcategory of $\Vect_{\C}^{\heart}$, all of its
objects are projective as well.
On the other hand,
$\Liq_p(\C)$ is
generated by compact projectives
(\cite[Theorem 6.5]{AnalyticCondensed}).
By \cite[Proposition 1.3.3.2]{HA} in conjunction with \cref{lem:
exact}, we therefore
obtain a canonical right t-exact functor 
of right-bounded derived
$\infty$-categories,
\[\un{\calD}: \calD^-(\Vect_{\C}^{f.d.}) \to \calD^-(\Liq_p(\C)).\]

\blem
\label{lem: fully faithful}
The functor $\un{\calD}$ is t-exact, fully faithful,
and symmetric monoidal.
\elem

\bproof
We first show that $\un{\calD}$ is fully faithful.
The categories $\calD^-(\Vect_{\C}^{f.d.})$ and $\calD^-(\Liq_p(\C))$
admit descriptions as $\rmN_{\dg}(\Ch^-(\Vect_{\C}^{f.d.}))$ and
$\rmN_{\dg}(\Ch^-({\Liq_p(\C)}_{\proj}))$, respectively,
where $\Ch^-(\calA)$ denotes the right-bounded chain
complexes in the additive category $\calA$ 
and $\rmN_{\dg}$ is the differential graded nerve of a differential
graded category\footnotemark.

\footnotetext{See \cite[Definition 1.3.1.1]{HA} for the definition of
differential graded category, \cite[Construction 1.3.1.6]{HA} for the
definition of differential graded nerve, and \cite[Definition 1.3.2.1]{HA} the
differential graded category structure on $\Ch^-(\calA)$.}

Let $V_n \in \Vect_{\C}^{f.d.} \subset \TVect^{\comp}_{\lc}$
be a vector space of dimension $n$, and choose a basis $\{v_i\}_{i \in S_n}$
for $V_n$, where $S_n$ is an index set of cardinality $n$. Then
$\un{V_n} \simeq \calM_{<p}(S_n)$, so $\un{V_n}$ is a
projective $p$-liquid vector space. This shows that
$\Vect_{\C}^{f.d.}$ actually embeds as a full subcategory
of ${\Liq_p(\C)}_{\proj}$. This embedding induces a
functor of differential graded categories,
\beqn
\label{eqn: map of Ch}
\Ch^-(\Vect_{\C}^{f.d.}) \to \Ch^-({\Liq_p(\C)}_{\proj}).
\eeqn
We claim that the functor 
given by the application
of $\rmN_{\dg}$ to 
\labelcref{eqn: map of Ch} is $\un{\calD}$.
To see this, note that by \cite[Theorem 1.3.3.2]{HA},
$\un{\calD}$ is the (essentially) unique extension $F$ of
$\un{(-)}: \Vect_{\C}^{f.d.} \to \Liq_p(\C)$ to derived
categories such that $\un{(-)} = 
\tau_{\geq 0} \circ (F|_{\calD^-(\Vect_{\C}^{f.d.})^{\heart}})$.
Now conclude by observing that 
$\rmN_{\dg}(\labelcref{eqn: map of Ch})$
satisfies this property.

Using \cite[Remark 1.3.3.6]{HA}, we
have that $\un{\calD}$ is t-exact because
the underlying functor on hearts, 
$\Vect_{\C}^{f.d.} \to \Liq_p(\C)$,
is exact by \cref{lem: exact}.
Finally, the symmetrical monoidality of
$\un{\calD}$ also follows immediately from
\cref{lem: exact}.
\eproof

\subsection{A theorem}

Before we state the theorem
of this section,
we state and prove
the following elementary lemma.

\blem
\label{lem: perf}
There is a t-exact equivalence of stable
$\infty$-categories,
\[\Perf \simeq \calD^b(\Vect_{\C}^{f.d.}).\]
\elem

\bproof
The inclusion $\Vect_{\C}^{f.d.}
\hook \Vect_{\C}^{\heart}$ induces by functoriality
a t-exact functor
\[\calD^-(\Vect_{\C}^{f.d.}) := \rmN_{\dg}(\Ch^-(\Vect_{\C}^{f.d.})) 
	\to \rmN_{\dg}(\Ch^-(\Vect_{\C}^{\heart})) =: \Vect^-_{\C}.\]
The restriction of this functor to the
subcategory of bounded complexes
induces a map, 
\[\calD^b(\Vect_{\C}^{f.d.}) \to \Perf\]
by the definition of $\Perf$. 
This functor is essentially
surjective by the definition of $\Perf$,
and is fully faithful because
the differential graded nerve
preserves fully faithfulness
(essentially because the Dold-Kan
functor preserves equivalences)
and the map of differential graded categories,
$\Ch^-(\Vect_{\C}^{f.d.}) \to 
\Ch^-(\Vect_{\C}^{\heart})$
is fully faithful. 
\eproof

We are now able to state and
prove our theorem.

\bthm
\label{thm: perfect functional calculus}
Suppose given $X \in \Perf$ and $T \in \End_{\Perf}(X)$.
Then for any entire function $f$, there exists an
endomorphism $f(T) \in \End_{\Perf}(X)$ such that
$H_i(f(T))$ is the linear map $f(H_i(T))$ obtained
from classical functional calculus by
applying $f$ to the induced map $H_i(T): H_i(X) \to H_i(X)$.
\ethm

\bproof
\Cref{lem: perf} in conjunction with \cref{lem: fully faithful}
tells us that $\Perf$ admits a t-exact fully
faithful embedding into $\calD^b(\Liq_p(\C))$.
On the other hand, the latter category is equivalent to
$\calD^b(\C, \calM_{<p})$ by \cref{prop: D Liq}, so we
have a fully faithful embedding,
$\Perf \hook \calD^b(\C, \calM_{<p})$
whose essential image clearly lies in
the full subcategory,
$\calD_{\Ban}(\C, \calM_{<p})$, 
of objects with homology objects belonging
to the essential image of $\Ban \hook \Liq_p(\C)$.
We denote this embedding by
\[\un{\calD}: \Perf \hook \calD_{\Ban}(\C, \calM_{<p}).\]
By the discussion in \cref{ssec: liquid functional calculus}
combined with \cref{prop: liquid-classical comparison},
$f(\un{\calD}(T))$ exists, and induces the endomorphisms
$f(H_*(\un{\calD}(T)))$ upon taking homology. Since
$\un{\calD}$ is fully faithful, $f(\un{\calD}(T))$ corresponds
to an endomorphism in $\Perf$ which we denote $f(T) \in \End_{\Perf}(X)$.
Finally, since $\un{\calD}$ is t-exact, $H_*(f(T))$ corresponds
under $\un{\calD}$ to the liquid endomorphism $f(H_*(\un{\calD}(T)))$.
\eproof

\printbibliography

\end{document}